\documentclass[12 pt,reqno]{amsart}
\usepackage{amstext, amssymb, amsthm, amsfonts, latexsym,bbm,mathtools}

\usepackage{enumerate}
\usepackage{caption} 
\usepackage[T1]{fontenc}
\usepackage[version=4]{mhchem}
\usepackage{ulem}
\usepackage{chemfig}
\usepackage{chemformula}
\usepackage{mathrsfs} 
\usepackage{MnSymbol}

\usepackage{comment}
\usepackage{tikz-cd}
\usepackage{subcaption}
\usepackage{floatrow}
\usepackage{hyperref}
\usepackage{geometry}

\usepackage{amsmath}
\let\emptyset\varnothing

\newtheorem{theorem}{Theorem}[section]
\newtheorem{lemma}[theorem]{Lemma}
\newtheorem{corollary}[theorem]{Corollary}
\newtheorem{proposition}[theorem]{Proposition}

\newtheorem{condition}{Condition}

\theoremstyle{definition}
\newtheorem{definition}[theorem]{Definition}
\newtheorem{example}[theorem]{Example}

\newcommand{\1}{\mathbf{1}}
\newcommand{\R}{\mathbb{R}}
\newcommand{\N}{\mathbb{N}}
\newcommand{\Z}{\mathbb{Z}}
\newcommand{\cC}{\mathcal{C}}
\newcommand{\cV}{\mathcal{V}}
\newcommand{\cR}{\mathcal{R}}

\newcommand{\cS}{\mathcal{S}}
\newcommand{\cX}{\mathcal{X}}

\newcommand{\cE}{\mathcal{E}}
\newcommand{\cN}{\mathcal{N}}

\newcommand{\cL}{\mathcal{L}}

\newcommand{\ess}{\mathrm{ess}}
\newcommand{\cyc}{\mathrm{cyc}}

\newcommand{\id}{\mathrm{id}}

\numberwithin{equation}{section}

\newcommand\red[1]{\textcolor{black}{#1}}

\begin{document}

\begin{figure}[h]
\href{https://doi.org/10.24072/pci.mcb.100405}{ 
\includegraphics[width=2in]{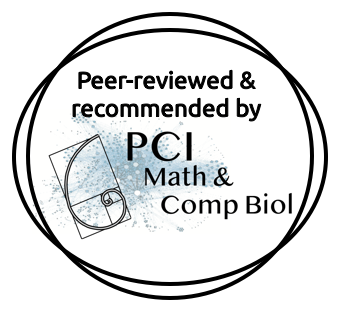}}
\end{figure}

\title[Reaction cleaving and complex-balanced distributions]{\red{Reaction cleaving and complex-balanced distributions for chemical reaction networks with general kinetics}}
\hfill\break

\author{
Linard Hoessly, Carsten Wiuf and Panqiu Xia}
\address{Data Center of the Swiss Transplant Cohort Study, University hospital Basel, Switzerland}
\email{\href{mailto:linard.hoessly@hotmail.com}{linard.hoessly@hotmail.com}}

\address{Department of Mathematical Sciences, University of Copenhagen, Denmark}
\email{\href{mailto:wiuf@math.ku.dk}{wiuf@math.ku.dk}}

\address{School of Mathematics, Cardiff University, UK}
\email{\href{mailto:xiap@cardiff.ac.uk}{xiap@cardiff.ac.uk}}
\date{}
\subjclass[2010]{}


\begin{abstract}
Reaction networks have become a major modelling framework in the biological sciences from epidemiology and population biology to genetics and cellular biology. In recent years, much progress has been made on stochastic reaction networks (SRNs),
modelled as continuous time Markov chains (CTMCs) and their stationary distributions. 
We are interested in complex-balanced stationary distributions, where the probability flow out of a  complex equals the flow into the complex.  We   characterise the existence and the form of complex-balanced distributions of SRNs with arbitrary transition functions through conditions on the cycles of the
 reaction graph (a digraph). Furthermore, we   give a sufficient condition for the existence of a complex-balanced distribution  and give precise conditions for when it is also necessary. The sufficient condition is also necessary for mass-action kinetics (and certain generalisations of that) or if the connected components of the digraph are cycles.  Moreover,  we state a deficiency theorem, a generalisation of the deficiency theorem for stochastic mass-action kinetics to arbitrary stochastic kinetics. The theorem gives the co-dimension of the parameter space for which a complex-balanced distribution exists.
 To achieve this, we construct an iterative procedure to decompose a strongly connected reaction graph into disjoint cycles, such that the corresponding SRN has equivalent dynamics and preserves complex-balancedness, provided the original SRN had so. This decomposition might have independent interest and might be applicable to  edge-labelled digraphs in general.
\end{abstract}


\maketitle

\section{Introduction}

Reaction networks offer a  framework to model  the dynamics of natural systems. 
They are applied across the sciences, for example in epidemiology \cite{murray,epid_CRN}, genetics \cite{ewens2004}, and   cellular biology \cite{wilkinson}.  A reaction network consists of a set of reactions, where a reaction represents a conversion, birth, or death of constituent particles (molecules, individuals, allele copies). For example, $A\ce{->} B$ might represent the conversion of one molecule of $A$ into one of $B$, and $S+I\ce{->} 2I$ might represent the infection of a susceptible individual by an infected individual, leading to two infected individuals.
  
A stochastic reaction network (SRN) is a homogeneous Markov chain on $\Z_{\ge0}^n$, given by an edge-labelled digraph (the reaction graph) \cite{AK15}, as illustrated in the example below:
 \begin{gather}\label{exp_sim}
\schemestart
  $A$
  \arrow{->[$\lambda_1$]}[30]
  $2C + D$
  \arrow{->[$\lambda_3$]}[-30]
  $B$,
  \arrow{<-[$\lambda_4$]}[210]
  $D$
  \arrow(@c1--@c3){<-[$\lambda_5$]}
  \arrow(@c1--@c4){->[$\lambda_2$]}
\schemestop
\end{gather}
The nodes (e.g., $A$) are  complexes, the edges represent reactions between complexes, and $\lambda_i\colon \Z_{\ge 0}^4\to\R_{\ge 0}$, $i=1,\ldots,4$, are the transition rates.  The vector of molecular counts of the species, $x=(x_A,x_B,x_C,x_D)\in\Z_{\ge 0}^4$, is the state of the system. If a reaction occurs, say, $A\ce{->} 2C + D$, then the Markov chain jumps from the current state $(x_A,x_B,x_C,x_D)$ to a new state $(x_A-1,x_B,x_C+2,x_D + 1)$; one molecule of $A$ is consumed, and two molecules of $C$ and one of $D$ are produced. 

The recent popularity of   
SRNs in the life sciences has lead to a deep interest in the existence and form of  stationary distributions  \cite{non-stand_1,anderson2,Cappelletti}. Analytical results are limited to  birth-death processes, finite state spaces, detailed- and complex-balanced systems with mass-action kinetics \cite{anderson2,kelly}, and some  special cases \cite{bibbona,Engblom09,dec_hoe}. In this paper, we discuss   complex-balanced stationary distributions with general kinetics (transition functions). Complex-balanced systems have their origin in Boltzman's work on detailed- (and cyclic) balanced systems, and have been the subject of much scrutiny      \cite{non-stand_1,anderson2,JoshiCap,Cappelletti,hong2023computational,hong2021derivation,kelly}:   a stationary distribution $\pi$ is complex-balanced if the probability flux out of a state through a   complex equals the flux into the state through the same complex, that is, if   
 \begin{align*}
  &\pi(x)\sum_{y'\colon \eta\to y' }\lambda_{\eta\to y'}(x)=\sum_{ y\colon y\to \eta }\pi(x+\phi(y)- \phi(\eta) ) \lambda_{y\to \eta}(x+\phi(y)- \phi(\eta) ),
  \end{align*}
 holds for all complexes    and   states   \cite{JoshiCap,Cappelletti}.    
 The function $\phi$   maps complexes to their stoichiometric coefficients, e.g., $\phi(A)=(1,0,0,0)$ and $\phi(2C + D)=(0,0,2,1)$, and the difference $\phi(y)- \phi(\eta)$ is the net molecular gain in the reaction $\eta\ce{->} y$.
   The existence of a complex-balanced stationary distribution implies the reaction graph is a disjoint union of strongly connected components (in the example, there is one such component) \cite{Cappelletti}. 

\red{
A main contribution of this paper is to construct a cleaving operation on reaction graphs that decomposes a strongly connected reaction graph of an SRN into a reaction graph of a dynamically equivalent SRN consisting of disjoint cycles only (a cyclic SRN). The construction is iterative, by splitting   nodes   with multiple incoming edges into multiple nodes with single incoming edges, and underlies  the proof of Theorem  \ref{thm_cleaved_RN} in particular. In order to formulate the procedure, we extend the definition of a `classical' SRN and allow multiple representations of the same complex and reaction (potentially with different labels) in the digraph.
Importantly, the cleaving operation preserves complex-balancedness:  A   cyclic  SRN is complex-balanced if and only if the original SRN is. All statements in the paper are valid for  classical SRNs as well as SRNs according to the new definition, since the classical definition is just a special case of the  new definition.}

\red{The cleaving operation might be of  independent interest. Furthermore, it is not restricted to reaction graphs of SRN, but  works on any strongly connected digraph with real labels, such as reaction graphs of deterministic reaction networks, with few changes. These changes relate to the discrete nature of the dynamics of SRNs versus the continuous nature of the dynamics of deterministic SRNs.}

In addition, we contribute the following:
\begin{itemize}
\item We characterise complex-balanced distributions of a reaction network with arbitrary kinetics through conditions on the cycles of its reaction graph
\item We provide a novel sufficient condition, extending   \cite[Theorem 4.1]{hong2023computational}, that implies existence of a complex-balanced distribution
\item Its necessity is also established when certain conditions are met, e.g. for mass-action kinetics and if the reaction graph is cyclic
\item We give   examples of SRNs for which we can find a stationary distribution by means of cycle decomposition
 \end{itemize}

Consider the SRN in \eqref{exp_sim}  that contains two cycles, $A\ce{->} 2C + D \ce{->} B \ce{->} A$ and $A\ce{->} D\ce{->} B \ce{->} A$. Splitting $A$ and $B$ into two nodes each, $(A,1)$, $(A,2)$, and $(B,1)$, $(B,2)$, respectively, results the following SRN,  
\begin{equation}\label{exp_sim1}
\schemestart
  $(A,1)$
  \arrow{->[$\lambda'_1$]}[30]
  $2C + D$
  \arrow{->[$\lambda'_3$]}[330]
  $(B,1)$
  \arrow(@c3--@c1){->[$\lambda'_5$]}
\schemestop
\quad \quad\quad
\schemestart
  $(A,2)$
  \arrow{->[$\lambda'_2$]}[30]
  $D$
  \arrow{->[$\lambda'_4$]}[330]
  $(B,2)$,
  \arrow(@c3--@c1){->[$\lambda'_6$]}
\schemestop
\end{equation}
where $(A,1)$ and $(A,2)$ are considered   different complexes  with the same stoichiometric coefficients, $\phi((A,1))=\phi((A,2))=(1,0,0,0)$, and likewise for $(B,1)$ and $(B,2)$. The $\lambda'_i$'s are kinetics to be defined, such that the SRN is dynamically equivalent to the original. In principle, this is not difficult, as one might take $\lambda'_i=\lambda_i$ for $i=1,\dots,4$, and choose arbitrary $\lambda_5'$ and $\lambda_6'$ with $\lambda_5'+\lambda_6'=\lambda_5$. However, this assignment does not necessarily preserve the complex-balanced property. If we choose $\lambda'_i=\lambda_i$ for $i=1,\dots,4$,
\begin{align*}
\lambda'_5(x)&=\frac{\lambda_1(x+e_A-e_B)\lambda_5(x)}{\lambda_1(x+e_A-e_B) + \lambda_2(x+e_A-e_B)}\1_{\{x'\colon x'_B\ge 1\}}(x),\\
\lambda'_6(x)&=\frac{\lambda_2(x+e_A-e_B)\lambda_5(x)}{\lambda_1(x+e_A-e_B) + \lambda_2(x+e_A-e_B)}\1_{\{x'\colon x'_B\ge 1\}}(x),
\end{align*}
where $e_A=(1,0,0,0)$ and $e_B=(0,1,0,0)$, then   the original SRN \eqref{exp_sim} is complex-balanced if and only if  the cleaved SRN \eqref{exp_sim1}   is   complex-balanced (see Section \ref{sec_pre} for the definition of cleaved SRNs).     Theorem \ref{thm_cleaved_RN} states the general procedure to decompose a strongly connected reaction graph into  disjoint cycles.

Graph decomposition techniques have   been used to study   deterministic reaction networks. Node balanced steady states  generalise complex-balanced steady states and are based on reaction graphs permitting multiple copies of the same complex \cite{feliu2018node}. Cyclic decompositions without dynamical equivalence have been constructed in \cite{Gopalkrishnan,horn1}. 
 A gluing operation was proposed in \cite{dec_hoe}.

The paper is organised as follows. In Section \ref{sec:main}, we discuss the main results of the paper; with applications given in Section \ref{sec_exam}.
In Section \ref{sec_pre}, we provide background on graphs and reaction networks and derive properties of decomposed reaction networks.  In Section \ref{sec_cleave}, building on the previous section, we introduce the cleaving operation that is used to decompose stochastic complex-balanced reaction networks into disjoint cycles. Finally,   proofs are in Section \ref{sec_prfs}.  

\section{Complex-balanced   stationary distributions}
\label{sec:main}

Let $(\cN,\lambda)$ be an SRN, where $\cN=(\cC,\cR)$ is a digraph of complexes and reactions, and $\lambda$ a labelling (kinetics) of the reactions (see Section \ref{sec_pre} for the precise definition), and let  $\Gamma\subseteq\Z_{\geq 0}^n$ be  a (closed) irreducible component of $(\cN,\lambda)$. We assume the following compatibility condition by default, which states that a reaction might `fire' only if the molecules of the source complex are available:

\begin{condition}\label{def_cptsto} 
 For $y\ce{->}y'\in \cR$ and $x\in\Gamma$, $\lambda_{y\to y'}(x)>0$ if any only if $x\geq \phi(y)$. %
\end{condition}

\begin{definition}
 A probability distribution $\pi$ on $\Gamma$ is a  
\begin{enumerate}[(i)]
\item {\bf stationary distribution} of $(\cN,\lambda)$, if for all $x\in \Gamma$,
\begin{align*}
&\pi(x)\sum_{y\to y'\in \cR}\lambda_{y\to y'}(x) =\sum_{y\to y'\in \cR}\pi\big(x-\phi(y')+\phi(y)\big)\lambda_{y\to y'}\big(x-\phi(y')+\phi(y)\big),
\end{align*}
where we set $\pi(x)=0$ and $\lambda_{y\to y'}(z)=0$ if $z\notin \Z_{\geq 0}^n$ (same below).

\item  {\bf complex-balanced distribution} of $(\cN,\lambda)$, if for all complexes $\eta\in \cC$, and all $x\in \Gamma$,
\begin{align}\label{def_vcbm}
&\pi(x)\sum_{y'\colon \eta\to y'\in \cR }\lambda_{\eta\to y'}(x)=\sum_{ y\colon y\to \eta\in \cR}\pi \big(x+\phi(y)-\phi(\eta)\big) \lambda_{y\to \eta}\big(x+\phi(y)-\phi(\eta)\big).
\end{align}

\item {\bf detailed-balanced distribution} of $(\cN,\lambda)$, if for all complexes $y,y'\in\cC$, and for  all $x\in \Gamma$,
\begin{align*}
\pi(x) \lambda_{y\to y'}(x) = \pi \big(x+\phi(y')-\phi(y)\big) \lambda_{y'\to y} \big(x+\phi(y')-\phi(y) \big),
\end{align*}
where $\lambda_{y\to y'}\equiv 0$ if $y \ce{->}y'\not\in\cR$.
\end{enumerate}
\end{definition}

A detailed-balanced distribution is also  complex-balanced,   and a complex-balanced distribution is also   stationary   \cite{JoshiCap}. Furthermore, complex-balancedness requires the digraph to be {\bf weakly reversible} (all connected components are strongly connected), that is,  every reaction  $y\ce{->}y'\in \cR$ belongs to a cycle $\gamma\subseteq \cR$ \cite{Cappelletti,craciun2020efficient}. Detailed balancedness requires the digraph to be {\bf reversible}, that is, if $y \ce{->}y' \in \cR$ then $y' \ce{->} y \in \cR$.

 Let $\cL_1,\dots, \cL_\ell$   be the   connected components of the digraph of  $(\cN, \lambda)$. Define
\begin{align}\label{def_gmmk}
\Gamma_k 
= \{x-\phi(y)\colon x\in \Gamma, y\in \cL_k\}\cap \Z_{\geq 0}^n, \quad k=1,\ldots,\ell.
\end{align} 
 
A main theorem is the following sufficient condition for the existence of a complex-balanced distribution.  

\begin{theorem}\label{prop_sdcrt}
Assume the digraph of  $(\cN, \lambda)$ is weakly reversible  with $\ell$ connected components.  Further, suppose   for each $k = 1, \dots, \ell$ and $y \to y' \in \cL_k$, that the kinetics factorises as   
\begin{align}\label{sdcrt}
\lambda_{y\to y'}(x)=\frac{\kappa_{y\to y'}}{m_k(x-\phi(y))g(x)},\quad x\geq \phi(y),\quad x\in\Gamma,
\end{align}
 where $g\colon\Gamma\to \R_{> 0}$ and $m_k\colon\Gamma_k\to \R_{> 0}$, $k=1,\dots,\ell$, are functions,  and $\kappa_{y\to y'}>0$ is a constant, such that for all complexes $\eta\in \cC$, 
\begin{align}\label{sdcrt0}
\sum_{ y'\colon \eta\to y'\in \cR }\kappa_{\eta\to y'}=\sum_{ y\colon y\to \eta\in \cR }\kappa_{y\to \eta}.
\end{align}
If $0<M \coloneqq \sum_{x\in\Gamma}g(x)<\infty$, then the distribution 
\begin{align*}
\pi(x)=\frac{1}{M}g(x),
\end{align*}
is  complex-balanced for $(\cN,\lambda)$.  
\end{theorem}

If the Markov chain is positive recurrent on $\Gamma$, then $\pi$ (and $g$ up to a normalising constant) is unique. Hence, also $m_k$ is unique up to a constant. If the Markov chain is transient and explosive, there might be several stationary distributions, hence the choice of $g$ and $m_k$ might not be unique.

Theorem \ref{prop_sdcrt}  
extends the condition presented in \cite[Theorem 4.1]{hong2023computational}, where  the $m_k$'s are required to be identical. The condition is also necessary is some cases:

\begin{proposition}\label{prop_cbmcr} 
Assume the digraph of  $(\cN, \lambda)$ is weakly reversible with   $\ell$ connected components.  Further, suppose that $\lambda_{y\to y'}(x)=\alpha_{y\to y'}\lambda_y(x)$ for $x\in\Gamma$ and $y\ce{->}y'\in \cR$, where $\alpha_{y\to y'} > 0$ is a constant. Then, a probability distribution $\pi$ on $\Gamma$ is  complex-balanced for $(\cN,\lambda)$, 
if and only if there exist non-negative functions $m_k$, $k=1,\dots, \ell$, such that 
\begin{align}\label{for_cpi}
\pi(x)=\kappa_{y\to y'}\big[\lambda_{y\to y'}(x)m_k(x-\phi(y))\big]^{-1},
\end{align}
for all $y\ce{->}y'\in \cR$, $y\in \cL_k$, and $x\in \Gamma$ with $x\geq \phi(y)$, where $\kappa_{y\to y'} > 0$  
satisfy   \eqref{sdcrt0}.
\end{proposition}

It seems surprisingly difficult to prove  Proposition \ref{prop_cbmcr}. We rely on Theorem \ref{thm_cleaved_RN} and the decomposition procedure developed in Section \ref{sec_cleave}. This procedure  provides a dynamically equivalent SRN with a digraph that consists of disjoint cycles only, while preserving the complex-balanced property.

In particular, Proposition \ref{prop_cbmcr} includes the case of stochastic mass-action kinetics,
\begin{equation*}
\lambda_{y\to y'}(x)=\alpha_{y\to y'}\frac{x!}{(x-y)!},\quad\alpha_{y\to y'}>0,\quad x\ge y,
\end{equation*}
and $0$ for $x\not\ge y$ (with $x!=\prod_{i=1}^nx_i$), where
\begin{equation}\label{eq:ma}
g(x)=\frac{c^x}{x!},\quad  m_k(x)=\frac{x!}{c^x},\quad \kappa_{y\to y'}= \alpha_{y\to y'}c^y,
\end{equation}
for some $c\in \R^n_{>0}$ (with $c^x=\prod_{i=1}^n c_i^{x_i}$). Also, the   kinetics proposed in \cite{non-stand_1,anderson3} fulfil the assumptions of Proposition \ref{prop_cbmcr}, and the form of the complex-balanced distributions could thus be found from the proposition. Other examples of kinetics fitting the framework of the assumptions of Proposition \ref{prop_cbmcr} are   stochastic Hill kinetics type I/II, and stochastic Michaelis-Menten kinetics \cite{anderson2,dec_hoe}. 

The correspondence between the two sets of constants in \eqref{eq:ma} and the requirement \eqref{sdcrt} are essential. For stochastic mass-action kinetics, the set of $\alpha_{y\to y'}$, $y\to y'\in\cR$, with this property has co-dimension   (also known as the deficiency)  $\delta=|\cC|-\ell -s$, where $|\cC|$ is the number of complexes,   
and $s$ the dimension of the space spanned by $\phi(y')-\phi(y)$, $y\to y'\in\cR$ \cite{Feinberg,craciun2009}.  The same holds for general kinetics (presented without proof):

\begin{lemma}\label{alternative_equivalence}
Assume the digraph of  $(\cN, \lambda)$ is weakly reversible  with $\ell$ connected components. Further,  suppose the kinetics factorises as   
\begin{align*}
\lambda_{y\to y'}(x)=\frac{\alpha_{y\to y'}}{m_k(x-\phi(y))g(x)},\quad x\geq \phi(y),\quad x\in\Gamma,
\end{align*}
  where  $g\colon\Gamma\to \R_{> 0}$ and $m_k\colon\Gamma_k\to \R_{> 0}$ are functions, and $\alpha_{y\to y'}>0$ is a constant,  such that there is  $c\in \R^n_{>0}$ fulfilling 
\begin{align*}
\sum_{ y'\colon \eta\to y'\in \cR }\alpha_{\eta\to y'}c^\eta=\sum_{ y\colon y\to \eta\in \cR }\alpha_{y\to \eta}c^y,\quad\text{for all}\quad \eta\in \cC.
\end{align*}
Then,   $\widehat g(x):=c^xg(x)$, $\widehat  m_k(x):=\frac{m_k(x)}{c^x}$ and   $\kappa_{y\to y'}:=\alpha_{y\to y'}c^y,$ $y\ce{->} y'\in\cR$ fulfil \eqref{sdcrt} and \eqref{sdcrt0}, and
 $$
\pi(x)=\frac{\widehat  g(x)}{\widehat  M}=\frac{c^xg(x)}{\widehat  M}, \quad x\in \Gamma,
$$ 
is  complex-balanced   for $(\cN, \lambda)$, provided $\widehat  M=\sum_{x\in\Gamma}\widehat g(x)<\infty$. The set of   $\alpha_{y\to y'}$, $y\ce{->} y'\in\cR$, for which this holds has co-dimension    $\delta=|\cC|-\ell -s$.
\end{lemma}

\begin{proposition}\label{prop_cbmpf}
Assume the digraph of  $(\cN, \lambda)$   consists of $\ell$ disjoint cycles. Then, a probability distribution $\pi$  on   $\Gamma$ is complex-balanced for $(\cN,\lambda)$, if and only if there exist  functions $m_k\colon\Gamma_k\to \R_{> 0}$, $k=1,\dots, \ell$,  such that for   $y\ce{->}y'\in \cL_k$, and   $x\in \Gamma$ with $x\geq \phi(y)$, we have
\begin{align}\label{for_pidcp}
\pi(x)=\big[\lambda_{y\to y'}(x)m_k(x- \phi(y))\big]^{-1}.
\end{align}
\end{proposition}

 We will show in Theorem \ref{thm_cleaved_RN} that any SRN with a weakly reversible digraph is complex-balanced if and only if it can be decomposed into an  SRN consisting of cycles only, such that Theorem \ref{prop_sdcrt} holds.   Thus, the cyclic SRN is also complex-balanced. In principle, one can therefore always use Proposition \ref{prop_cbmpf} to determine the stationary distribution of a complex-balanced SRN by this cyclic decomposition.

In general, one might be able to decompose the digraph of an SRN into disjoint cycles while preserving dynamical equivalence in many ways, as already alluded to in example \eqref{exp_sim1}. However, the difficulty does not lie in preserving dynamical equivalence, but in preserving the complex-balance property. Theorem \ref{thm_cleaved_RN} provides one way of achieving this.

\begin{example}\label{exp_trag} 
 Consider the mass-action SRN,
\begin{equation}\label{cexp_21}
\schemestart 
  $C$,
  \arrow{<=>[$2$][$1$]}[150]
  $A$
  \arrow{<=>[$2$][$1$]}[210]
  $B$
  \arrow(@c3--@c1){<=>[$1$][$2$]}
\schemestop
\end{equation}
and let $\Gamma$ be an irreducible component. Then, $\pi(x)=\frac{M_{\Gamma}}{x!}$, where  $M_{\Gamma} > 0$ a normalisation constant, is the unique complex-balanced distribution for \eqref{cexp_21}  
\cite[Theorem 4.1]{anderson2}.
A dynamically equivalent   SRN, decomposed according to Theorem \ref{thm_cleaved_RN}, with five disjoint cycles and mass-action kinetics is (details omitted):
\begin{equation*}
\schemestart
  $(A,1)$
  \arrow{<=>[$\alpha_1$][$\alpha_2$]}
  $(B,1)$, 
\schemestop
\qquad
\schemestart
  $(B,2)$
  \arrow{<=>[$\alpha_3$][$\alpha_4$]}
  $(C,2)$,
\schemestop
\qquad
\schemestart
$(A,3)$
  \arrow{<=>[$\alpha_5$][$\alpha_6 $]}  
  $(C,3)$,
  \schemestop
\end{equation*}
\begin{equation*}
\schemestart
  $(C,4)$
  \arrow{->[$\alpha_7$]}[150]
  $(A,4)$
  \arrow{->[$\alpha_8$]}[210]
  $(B,4)$
  \arrow(@c3--@c1){->[$\alpha_9$]}
\schemestop
\quad \mathrm{and}\quad
\schemestart
   $(C,5)$,
  \arrow{<-[$\alpha_{10}$]}[150]
  $(A,5)$
  \arrow{<-[$\alpha_{11}$]}[210]
  $(B,5)$
  \arrow(@c3--@c1){<-[$\alpha_{12}$]}
\schemestop
\end{equation*}
where $\alpha_1, \dots,\alpha_{12}$ are rate constants satisfying $\alpha_1=\dots=\alpha_6=1-\beta$, $\alpha_7=\alpha_8=\alpha_9=\beta$, and $\alpha_{10}=\alpha_{11}=\alpha_{12}=\beta+1$, with  $\beta\in(0,1)$ arbitrary. Then, the cyclic SRN is complex-balanced. The   procedure in Section \ref{sec_cleave} results in $\beta=\frac{1}{15}$. 
\end{example} 

\begin{example}\label{cexp}
Not all decompositions into cycles preserve the property of being complex-balanced. Consider the mass-action SRN,
\begin{equation}\label{cexp_1}
A\ce{<=>[3][3]}B.
\end{equation}
  The following system with   mass-action kinetics
\begin{equation}\label{cexp_2}
(A,1)\ce{<=>[1][2]}(B,1),\quad (A,2)\ce{<=>[2][1]}(B,2),
\end{equation} 
and $\phi ((A,i))=A$ and $\phi((B,i))=B$ for $i=1,2$, is a cleaved SRN of \eqref{cexp_1}. The probability distribution $\pi(x)=\frac{M_{\Gamma}}{x!}$, where $M_\Gamma$ is a normalising constant, is a complex-balanced distribution for \eqref{cexp_1} on any irreducible component $\Gamma$ \cite[Theorem 4.1]{anderson2}; hence also a stationary distribution. But $\pi$ is not a complex-balanced distribution of \eqref{cexp_2}; only a stationary distribution.
\end{example}

 We end the section with a version of  Theorem \ref{prop_sdcrt}  for detailed-balanced SRNs. It is sufficient to consider only cycles of length two:

\begin{proposition}\label{prop_DB}
Assume the digraph of  $(\cN, \lambda)$ is   reversible. Then, a probability distribution $\pi$ on $\Gamma$ is detailed-balanced 
if and only if there exist functions $m_{y \to y'} \colon \Gamma_{y,y'} \to \R_{>0}$ for all $y \ce{<=>} y' \in \cR$, such that $m_{y \to y'} = m_{y' \to y}$ and
\begin{align}\label{eq_prop_db}
\lambda_{y \to y'}(x) = [m_{y \to y'}(x - \phi(y)) \pi(x)]^{-1},
\end{align}
where $\Gamma_{y,y'}  = (\{x-\phi(y)\colon x\in \Gamma\}\cup\{x-\phi(y')\colon x\in \Gamma\})\cap \Z_{\geq 0}^n.$
\end{proposition}

\section{Examples}\label{sec_exam}
 
    We present some examples to illustrate the results.  Most reaction networks used in applications in biophysics, cellular biology and systems biology are not weakly reversible, let alone reversible. To remedy this, various techniques of network translation have been invented, that is, ways to transform a non-weakly reversible reaction network into a dynamically equivalent weakly reversible reaction network \cite{hong2023computational,hong2021derivation,johnston1,johnston}. For example, one might add or delete  species in equal numbers on both sides of a reaction, or split a reaction into two while preserving the total reaction rate. We will make use of these techniques too.  
    
The main aim is to illustrate Theorem \ref{prop_sdcrt}.   The one-node cleaving procedure used to prove Proposition \ref{prop_cbmcr}, is generally very laborious to apply (as there might be many cycles and nodes are cleaved one by one). We give one example of the procedure. In many cases, it seems more appropriate to adopt a manual approach to cleaving.  

\subsection{Michaelis-Menten kinetics} Consider an enzyme-regulated mechanism for product formation with Michaelis-Menten kinetics   \cite{cornish-bowden}: 
\begin{equation*}
S+E^*\ce{<=>[\lambda_1][\lambda_2]} P+E, \quad 0\ce{<=>[\lambda_3][\lambda_4]} E.
\end{equation*}
A substrate $S$ is converted into a product $P$ reversibly by means of an active enzyme $E^*$, which is then converted to its `inactive' form $E$. The enzyme might further be supplied from the surroundings and degraded. The reactions happen in a density-regulated manner:
\begin{align*}
\lambda_1(x)&=\frac{\alpha_1x_Sx_{E^*}}{\beta_1+x_S},\quad\lambda_2(x) =\frac{\alpha_2x_Px_E}{\beta_2+x_P},\quad \lambda_3(x) =\frac{\alpha_3}{\beta_3+x_E},\quad\lambda_4(x) =\frac{\alpha_4x_E}{\beta_4+x_E},
\end{align*}
where $\alpha_i, \beta_i$, $i=1,\ldots,4$, are positive constants and $\Gamma_T=\{ x\in\Z^4_{\ge0}\colon x_S + x_P=T, x_E\ge 0, x_{E^*}\ge0\}$ is an irreducible component for $T\in\N$. All of the statements are applicable in this case. We   apply     Lemma \ref{alternative_equivalence} to achieve
\begin{align*}
m_1(x)&=\frac{x!}{ \prod_{i=0}^{x_S}(\beta_1+i)\prod_{i=0}^{x_P}(\beta_2+i)},\\
m_2(x)&=\frac{x!}{(\beta_3+x_E)\prod_{i=0}^{x_S}(\beta_1+i)\prod_{i=0}^{x_P}(\beta_2+i)},\\
g(x)&= \frac{ \prod_{i=0}^{x_S}(\beta_1+i)\prod_{i=0}^{x_P}(\beta_2+i)}{x!},
\end{align*}
provided $\beta_3=\beta_4$. Since, there exists $c\in\R^4_{>0}$, such that $\alpha_1c_Sc_{E^*}=\alpha_2c_Pc_E$ and $\alpha_3=\alpha_4c_E$, then the conditions of Lemma \ref{alternative_equivalence}  are fulfilled, and $c^x g(x)$ normalised is a detailed-balanced distribution for all $\alpha_i>0$, $i=1,\ldots,4$, and $\beta_1,\beta_2>0$, provided $\beta_3=\beta_4>$. Hence, the co-dimension for the $\alpha_i$-parameters is $\delta=|\cC|-\ell-s=4-2-2=0$.

\subsection{Phosphorylation mechanism} 

Consider  a  phosphorylation mechanism  modelled  with mass-action kinetics \cite{hong2021derivation}:
\begin{equation*}\label{pak1}
A\ce{<-[\alpha_1]} B\ce{<=>[\alpha_2][\alpha_3]} C,\quad 2A\ce{->[\alpha_4]} A+B,\quad A+B\ce{<=>[\alpha_5][\alpha_6]} A+C,
\end{equation*}
  on the irreducible component $\Gamma_T=\{(x_A,x_B,x_C)\in \Z_{\ge 0}^3 \colon x_A+ x_B+x_C=T,x_A\ge 1\}\subseteq \Z_{\ge 0}^3$. Here, $A$ is a substrate with two phosphorylation sites free, $B$ has one of these sites occupied by a phosphate group, and $C$ has both sites occupied. The different reactions represent different ways phosphorylation occurs. 
  
  The SRN is not weakly reversible, so we modify it into a dynamically equivalent weakly reversible SRN by adding $A$ on both sides in the first three reactions, and decomposing the SRN according to Theorem \ref{thm_cleaved_RN}. Then, this dynamically equivalent weakly reversible SRN is complex-balanced if and only if the below SRN is:
\begin{equation*}\label{pak1-2}
2A\ce{<=>[\lambda_4][\lambda_1]}(1, A+B)\ce{<=>[\lambda_2][\lambda_3]}(1,A+C),\quad (2,A+B)\ce{<=>[\lambda_5][\lambda_6]} (2,A+C)
\end{equation*}
$$\lambda_1(x)  =\alpha_1x_B,\quad \lambda_2(x) =\alpha_2 x_B,\quad \lambda_3(x) =\alpha_3 x_C,$$
 $$\lambda_4(x)  =\alpha_4x_A (x_A-1),\quad \lambda_5(x)  =\alpha_5 x_Ax_B,\quad \lambda_6(x) =\alpha_6 x_Ax_C.$$
 for $x\in\Gamma$. Condition \ref{def_cptsto} is fulfilled on $\Gamma_T$. Choosing
$$m_1(x)=\frac{(x_A+1)!x_A! x_B!x_C!}{(\alpha_1\alpha_3)^{x_A}(\alpha_3\alpha_4)^{x_B}(\alpha_2\alpha_4)^{x_C}},\quad m_2(x)=\frac{(x_A!)^2 x_B!x_C!}{(\alpha_1\alpha_3)^{x_A}(\alpha_3\alpha_4)^{x_B}(\alpha_2\alpha_4)^{x_C}},$$
$$\kappa_1=\alpha_1^2\alpha_3^2\alpha_4,\quad \kappa_2=\alpha_1 \alpha_2\alpha_3^2\alpha_4,\quad \kappa_3=\alpha_1 \alpha_2\alpha_3^2\alpha_4,$$
$$ \kappa_4=\alpha_1^2\alpha_3^2\alpha_4,\quad\kappa_5=\alpha_1\alpha_3^2\alpha_4\alpha_5,\quad \kappa_6=\alpha_1\alpha_2\alpha_3\alpha_4\alpha_6,$$
ensures the two components fulfil  Theorem \ref{prop_sdcrt} with  complex-balanced distribution:
$$\pi_T(x)=M_T\frac{(\alpha_1\alpha_3)^{x_A}(\alpha_3\alpha_4)^{x_B}(\alpha_2\alpha_4)^{x_C}}{ x_A!(x_A-1)!x_B!x_C!},\quad x\in\Gamma_T,$$
provided $\kappa_5=\kappa_6$, that is, if $\alpha_3\alpha_5=\alpha_2\alpha_6$; 
where $ M_T>0$ is a constant. 

The second component has the mass-action form, with the number of $A$'s being conserved. Hence, the stationary distribution should take the form in \eqref{eq:ma}. It also does so: treating $x_A$ as constant, then $\pi_T$ has mass-action form.

We note that the factorisation of the kinetics takes the form in Lemma \ref{alternative_equivalence}. In particular, the co-dimension is $\delta=|\cC|-\ell-s=5-2-2=1$, hence there is one constraint on the parameters giving rise to complex-balanced distributions, as also found.

\subsection{Modified birth-death process}

The following example is a modification of a classical birth-death process that has an extra reaction with a jump of size two \cite{ACGW15}. More precisely, we consider the following SRN on $\Gamma=\N_0$ with mass-action kinetics, which is not weakly reversible,
\begin{equation}\label{exp_bd}
A \ce{<=>[\alpha_1][\alpha_2]} 0 \ce{->[\alpha_3]} 2A.
\end{equation}

To apply Theorem \ref{prop_sdcrt} we need to find an equivalent weakly reversible SRN.
 Changing the reaction $A\ce{->} 0$ to $A\ce{->}0$ and $2A\ce{->}A$, then we look for a dynamically equivalent SRN of the following form,
\begin{equation}\label{exp_bdt}
\schemestart
  $A$
  \arrow{<=>[$\lambda_1$][$\lambda_2$]}[30]  
  $0$
  \arrow{->[$\lambda_3$]}[-30]
  $2A$
  \arrow(@c1--@c3){<-[$\lambda_4$]}
\schemestop
\end{equation}
\begin{align*}
\lambda_2 (x) = \alpha_2, \quad \lambda_3 (x) = \alpha_3,\quad 
 \lambda_4 (x) + \lambda_1 (x) = \alpha_1 x,
 \end{align*}
 where $x$ denotes the number of $A$ molecules. We will show that $\lambda_1$ and $\lambda_4$ are uniquely determined for any fixed $\alpha_1,\alpha_2,\alpha_3\in \R_{>0}$ such that \eqref{exp_bdt} is complex-balanced and Condition \ref{def_cptsto} is satisfied.  
 
The SRN \eqref{exp_bdt} fulfils the `constant ratio' condition in Proposition \ref{prop_cbmcr}, hence Theorem \ref{prop_sdcrt} can be applied to justify complex-balancedness. However, we prefer to decompose it into cycles to avoid the difficulty of choosing the constants $\kappa$'s in \eqref{sdcrt0} when applying Proposition \ref{prop_cbmpf}. 
Cleaving the SRN into cycles results in
\begin{equation}\label{exp_bdtsplt}
\begin{aligned}[c]
\cL_1:&\quad (A, 1) \ce{<=>[\lambda_{\cyc,1}][\lambda_{\cyc,2}]} (0, 1),
\end{aligned}\quad \qquad
\begin{aligned}
\cL_2:
\end{aligned}\quad 
\begin{aligned}
\schemestart
  $(A, 2)$
  \arrow{->[$\lambda_{\cyc,5}$]}[30]  
  $(0, 2)$
  \arrow{->[$\lambda_{\cyc,3}$]}[-30]
  $(2A, 2)$,
  \arrow(@c1--@c3){<-[$\lambda_{\cyc,4}$]}
\schemestop
\end{aligned}
\end{equation}
where $\lambda_{\cyc, i} = \lambda_{i}$ for $i = 2,3,4$,
\begin{align*}
\lambda_{\cyc,1} (x) = \frac{\alpha_2}{\alpha_2 + \alpha_3}\lambda_1 (x), \quad \mathrm{and}\quad \lambda_{\cyc,5} (x) = \frac{\alpha_3}{\alpha_2 + \alpha_3}\lambda_1 (x).
\end{align*}
Due to Proposition \ref{prop_cbmpf} and Theorem \ref{thm_cleaved_RN}, the SRN \eqref{exp_bdt} is complex-balanced, if and only if there exist non-negative functions $m_1$, $m_2$ and $g$ on $\Z_{\geq 0}$, such that
\begin{align}\label{eq_bdt1}
\lambda_{\cyc, 1} (x + 1) g(x + 1) = \lambda_{\cyc,2} (x) g(x) = m_1 (x)^{-1},
\end{align}
and
\begin{align}\label{eq_bdt2}
\lambda_{\cyc, 3} (x) g(x) = \lambda_{\cyc,4} (x + 2) g(x + 2) = \lambda_{\cyc, 5} (x + 1) g(x + 1) = m_2 (x)^{-1}
\end{align}
for all $x\in \Z_{\geq 0}$ and $M\coloneqq \sum_{x = 1}^{\infty} g(x) \in (0,\infty)$. If we choose 
\begin{align}\label{eq_bdt5}
m_1 (x) = \alpha_3 m_2 (x)/\alpha_2,
\end{align}
then \eqref{eq_bdt1} is a consequence of \eqref{eq_bdt2}, and we only need to solve for \eqref{eq_bdt2}. Suppose \eqref{eq_bdt2} holds. For all $x\geq 1$, it follows that
\begin{align*}
\frac{\lambda_{\cyc, 3}(x)}{\lambda_{\cyc, 5} (x+1)} = \frac{\lambda_{\cyc, 5} (x)}{\lambda_{\cyc, 4} (x+1)},
\end{align*}
 and thus,
\[
\lambda_1 (x + 1) = \frac{\alpha_1 (\alpha_2 + \alpha_3)^2 (x + 1)}{(\alpha_2 + \alpha_3)^2 + \alpha_3 \lambda_1 (x)}.
\]
Condition \ref{def_cptsto} gives $\lambda_{4} (0) = \lambda_{4} (1) = 0$. Thus, $\lambda_{1} (1) = \alpha_1$ and $\lambda_1$ is uniquely determined by the recursion:
\begin{align}\label{for_rec}
\lambda_1 (x) = 
\begin{cases}
\displaystyle 0, & x = 0,\\
\displaystyle \frac{\alpha_1 (\alpha_2 + \alpha_3)^2 x}{(\alpha_2 + \alpha_3)^2 + \alpha_3 \lambda_1 (x-1)}, & x > 0.
\end{cases}
\end{align}
In fact, in \eqref{for_rec}, $0 < \lambda_1 (x) < \alpha_1 x$ whenever $\lambda_1 (x-1) > 0$. Therefore, $\lambda_1 (x)$ and $\lambda_4(x)=\alpha_1 x-\lambda_1(x)$ are in $(0,\alpha_1 x)$ for all $x\geq 2$, and Condition \ref{def_cptsto} holds for $\lambda_{1}$ and $\lambda_4$. Assume $g(0) = 1$, then combined with \eqref{eq_bdt2}, we have
\begin{align}\label{eq_bdt6}
g (x + 1) = \frac{\lambda_{\cyc,3} (x)}{\lambda_{\cyc, 5} (x + 1)} g(x) = \frac{\alpha_2 + \alpha_3}{\lambda_1 (x + 1)} g(x) = (\alpha_2 + \alpha_3)^{x + 1} \Big(\prod_{u= 1}^{x + 1} \lambda_1 (u)\Big)^{-1}
\end{align}
and
\begin{align}\label{eq_bdt4}
m_2 (x) = (\alpha_3 g (x))^{-1},
\end{align}
for all $x\geq 0$. With $\lambda_1$, $m_1$ $m_2$ and $g$ defined as in \eqref{for_rec}, \eqref{eq_bdt5}, \eqref{eq_bdt4} and \eqref{eq_bdt6}, respectively, one can verify that Theorem \ref{prop_sdcrt}(ii) is satisfied. 

To prove the existence of a complex-balanced distribution, we need to show $M \coloneqq \sum_{x = 1}^{\infty} g(x) <\infty$. 
By using the recursive formula \eqref{for_rec}, we deduce that for all $x \geq 1$,
\[
\lambda_1 (x)\lambda_1 (x + 1) = (x + 1) h \big( \lambda_1 (x)\big), \quad \mathrm{where}\quad h (u)\coloneqq \frac{\alpha_1 (\alpha_2 + \alpha_3)^2 u}{(\alpha_2 + \alpha_3)^2 + \alpha_3 u}, \quad u\in \R_{\geq 0}.
\]
Due to \eqref{for_rec} and the fact that $0\leq \lambda_1 (x) \leq \alpha_1 x$, we have for $x \geq 2$,
\[
\lambda_1 (x) \geq \frac{\alpha_1 (\alpha_2 + \alpha_3)^2 x}{(\alpha_2 + \alpha_3)^2 + \alpha_3 \alpha_1 (x-1)} \geq \frac{\alpha_1 (\alpha_2 + \alpha_3)^2 (x - 1)}{(\alpha_2 + \alpha_3)^2 + \alpha_3 \alpha_1 (x-1)} \geq c_0 \coloneqq \frac{\alpha_1 (\alpha_2 + \alpha_3)^2}{(\alpha_2 + \alpha_3)^2 + \alpha_3 \alpha_1},
\]
where the last inequality follows from the property that $x\mapsto \frac{ax}{c+bx}$ is increasing on $\R_{\geq 0}$ with arbitrary parameters $a,b,c>0$. Since $h$ is also increasing on $\R_{\geq 0}$, it holds that for all $x\geq 2$,
\begin{align*}
\lambda_1 (x)\lambda_1 (x + 1) \geq (x + 1) \inf_{x\geq 2} h \big(\lambda_1 (x)\big) \geq h(c_0)>0.
\end{align*}
As a consequence, for $x \geq \frac{2 (\alpha_2 + \alpha_3)^2}{h(c_0)} \vee 2$,
\[
\frac{g (x+1)}{g(x-1)} = \frac{(\alpha_2 + \alpha_3)^2}{\lambda_1 (x+1)\lambda_1(x)} \leq \frac{(\alpha_2 + \alpha_3)^2}{(x + 1) h(c_0)} < \frac{1}{2},
\]
and $M$ is finite by the ratio test. Due to Proposition \ref{prop_cbmpf} and Theorem \ref{thm_cleaved_RN}, $\pi(x) \coloneqq \frac{1}{M} g(x)$ is the unique complex-balanced distribution for \eqref{exp_bdtsplt} and also \eqref{exp_bdt}, and thus a stationary distribution for \eqref{exp_bd}. From \cite{xu}, the reaction network \eqref{exp_bd} is positive recurrent, hence this distribution is the unique stationary distribution.

\section{Stochastic reaction networks}\label{sec_pre}

In this section, we define SRNs and present a partial order on the space of SRNs. The main decomposition theorem (cleaving of SRNs) will make use of this partial order.

\subsection{Notation} 

Let $\R$, $\R_{\geq 0}$ and $\R_{>0}$ be the set of real, non-negative and positive numbers, respectively. Let $\Z$ and $\Z_{\geq 0}$ 
 be the set of integers and non-negative integers, respectively. 
For $x=(x_1,\dots, x_n), y=(y_1,\dots, y_n)\in \R^n$, we define $x\geq y$, if $x_i\geq y_i$ for all $i=1,\dots, n$; and $x>y$ if $x\geq y$ and $x\neq y$.   Furthermore, for $x\in \R_{\ge 0}^n$, $y\in \Z^n_{\ge 0}$, the notation $x^y$ is used for $\prod_{i=1}^nx_i^{y_i}$, and for $x\in \Z_{\geq 0}^n$, we write $x!$ for $x_1!\cdots x_n!$.

\subsection{Graph theory}

Consider a {\bf digraph}  $(\cV,\cE)$, where $\cV$ is a finite set of {\bf nodes} and $\cE\subseteq \cV\times \cV$ is a finite set of {\bf edges}. A {\bf sub-digraph} $(\cV',\cE')$ of $(\cV,\cE)$ is a digraph such that $\cV'\subseteq \cV$ and $\cE'\subseteq (\cV'\times \cV')\cap \cE$. Two sub-digraphs are {\bf disjoint} if their sets of nodes are disjoint. 

A {\bf walk} is an ordered finite sequence of edges in $(\cV,\cE)$, denoted by $\theta=(v_1\ce{->}v_2, v_2\ce{->}v_3,\dots, v_{k-1}\ce{->}v_k)$ or $(v_1\ce{->}v_2\ce{->}\cdots\ce{->}v_k)$ for convenience. The walk is {\bf closed} if $v_1=v_k$, and is {\bf open} if it is not closed. An open walk $\theta$ is {\bf directed} from $v_1$ to $v_k$, and {\bf links} $v_1$ and $v_k$, vice versa $v_k$ and $v_1$. 
If all nodes are different, then $\theta$ is a {\bf path}, and if all nodes are different but $v_1=v_k$, then it is a {\bf cycle}.  Paths and cycles, but not walks, might be seen as sub-digraphs.

 A digraph $(\cV,\cE)$ is {\bf connected}, if for any pair of nodes $v,v'\in \cV$, there exist nodes $v_0,v_1,\dots, v_{k+1}$, and paths $\theta_1,\dots, \theta_{k+1}$ in $(\cV, \cE)$, such that $v_0=v'$ and $v_{k+1}=v$ and $\theta_i$ links $v_{i-1}$ and $v_i$ for all $i=1,\dots, k+1$. A  sub-digraph $(\cV',\cE')$ is a {\bf connected component} of $(\cV,\cE)$, if $\cE'= (\cV'\times \cV')\cap \cE$ and no nodes $v\in\cV\setminus\cV'$ are linked to a node in $\cV'$.  A connected component is {\bf strongly connected} if there is  a path from $v$ to $v'$ for any pair of nodes $v,v'\in  \cV'$.

 The ensuing lemma then follows by definition.
\begin{lemma}\label{lmm_cycle} 
Let $(\cV,\cE)$ be a digraph satisfying the following
\begin{enumerate}[(i)]
\setlength\itemsep{0.5em}
\item For any edge $e\in \cE$ there is a cycle $\gamma\subseteq \cE$ with $e\in\gamma$.

\item For any node $v\in \cV$ there is at most one edge $e\in \cE$ such that $e=w\ce{->} v$ with $w\in\cV$. 
\end{enumerate} 
Then, $(\cV,\cE)$ consists of disjoint cycles. 
\end{lemma}

\subsection{SRNs} 

In our context, an SRN is a pair $(\cN,\lambda)$, where $\cN$ is a  digraph $(\cC,\cR)$ on a set $\cS=\{S_1,\ldots,S_n\}$  with a map 
\[
\phi\colon\cC\to\Z_{\geq 0}^{\cS},\quad\Z_{\geq 0}^{\cS}=\Big\{\sum_{i=1}^nz_i S_i\colon z_i\in \Z_{\geq 0},  i=1,\ldots,n\Big\}\cong \Z^n_{\ge0},
\]
that associates to each node a non-negative integer vector.
The elements of $\cS$ are {\bf species}, those of $\cC$ are {\bf complexes}, and those of $ \cR$ are {\bf reactions}.
For a reaction $r=y\ce{->}y'\in \cR$, $y,y'\in\cC$, the node $y$ is 
 the {\bf reactant} and $y'$ the {\bf product}. Moreover, $r$ is called an {\bf incoming} reaction of complex $y'$, and an {\bf outgoing} reaction of complex $y$. The vector $\phi(y)$ gives the {\bf species composition} of  a complex $y$.

Furthermore, \[\lambda=(\lambda_{y\to y'} \colon y\ce{->}y' \in \cR),\quad \lambda_{y\to y'}\colon\Z_{\geq 0}^n\to \R_{\geq 0},
\]
is an {\bf edge-labelling} of the digraph, referred to as the {\bf   kinetics}. The evolution of the species counts $X(t)$, $t\ge0$, over time is modelled as a $\Z_{\geq 0}^n$-valued CTMC, satisfying the following SDE:
\begin{align}\label{def_dynsto}
X(t)=X(0)+\sum_{y\to y'\in \cR}Y_{y\to y'}\Big(\int_0^t\lambda_{y\to y'}(X(s))ds\Big)(\phi(y')-\phi(y)),
\end{align} 
where $Y_{y\to y'}$, $ y\ce{->}y'\in \cR$, is a collection of i.i.d. unit rate Poisson processes; that is, $\lambda_{y\to y'}$ is the transition intensity at which reaction $y\to y'$ `fires'. 
 When writing $(\cN,\lambda)$, we implicitly assume $\cC,\cR,\cS,\phi$ are given. On occasion,  we write  for brevity, $\cN=(\cC,\cR,\cS,\phi)$.

The two graphs \eqref{exp_sim} and \eqref{exp_sim1} are  SRNs. In the first, $\phi=\id_{\cC}$, while in the second $\phi'\colon\cC\to \Z_{\geq 0}^4$ is given by $\phi'((y,j))=y$ for $y\in\{A,D\}$, $j=1,2$, and $\phi'(y)=y$ for $y\in\{B,C\}$.

\begin{definition}
Let $(\cN,\lambda)$ and $(\cN',\lambda')$ be two SRNs. If there exists a map $\psi\colon\cC'\to\cC$, such that
\begin{enumerate}
\item[(i)] $\phi'=\phi\circ \psi$, which implies  $\cC=\psi(\cC')=\{\psi(y)\colon y\in \cC'\}$. 

\item[(ii)] $\cR=\psi(\cR')=\{\psi(y)\ce{->}\psi(y')\colon y\ce{->}y'\in\cR'\}$.
\item[(iii)] For all $x\in \Z_{\geq 0}^n$ and all $y\ce{->}y'\in \cR$,  
\begin{align}\label{equ_kcleaved SRN}
\lambda_{y\to y'}(x)=\sum_{r\in \psi^{-1}(y\to y')}\lambda'_{r}(x).
\end{align}
\end{enumerate}
Then, $(\cN',\lambda')$ is a {\bf cleaved SRN}  of $(\cN,\lambda)$ with {\bf projection} $\psi$, denoted  $(\cN', \lambda') \succeq (\cN, \lambda)$. 
 If only (i)-(ii) hold, then we write $\cN'\succeq\cN$. Furthermore,   complexes  $y'\in \cC'$ for which $\psi(y')=y\in \cC$ are called  {\bf copies} of $y$. 
 \end{definition}
 
Digraph \eqref{exp_sim1} is  a cleaved SRN of digraph \eqref{exp_sim} with projection $\psi=\phi'$, provided (iii) is fulfilled. 
 `Being cleaved' is a partial order on the set of SRNs. 

\begin{lemma}
Let $(\cN,\lambda)$, $(\cN',\lambda')$ and $(\cN'',\lambda'')$ be SRNs. Suppose  $(\cN',\lambda')\succeq (\cN,\lambda)$ with  projection $\psi$, and $(\cN'',\lambda'')\succeq (\cN',\lambda')$ with the projection $\psi'$. Then, $(\cN'',\lambda'')\succeq (\cN,\lambda)$ with  projection $\psi\circ\psi'$.
\end{lemma}

 The {\bf essential SRN}  $(\cN_\ess,\lambda_{\ess})$ of an SRN $(\cN,\lambda)$ is defined by  $\cN_{\ess} \coloneqq  (\phi(\cC),\phi(\cR),\cS,\text{id}_{\phi(\cC)})$, where
\begin{align*}
\phi(\cC) \coloneqq\{\phi(y)\colon y\in \cC\},\quad
\phi(\cR) \coloneqq \{\phi(y)\ce{->}\phi(y')\colon y\ce{->} y'\in \cR\},
\end{align*}
and $\text{id}_{\phi(\cC)}$ is the identity map on $\phi(\cC)$; and 
 \[
\lambda_{\ess, y\to y'}(x) \coloneqq \sum_{r\in \phi^{-1}(y\to y')}\lambda_{r}(x)
\]
for all $y\ce{->} y'\in \phi(\cR)$ and $x\in \Z_{\geq 0}^n$. Clearly, $(\cN,\lambda) \succeq (\cN_{\ess}, \lambda_{\ess})$ with projection $\psi=\phi$, and $((\cN_\ess)_\ess, (\lambda_{\ess})_{\ess}) = (\cN_{\ess}, \lambda_{\ess})$.

\begin{lemma}\label{lmm_drsplt}
Suppose $(\cN',\lambda')\succeq (\cN,\lambda)$. Then, $(\cN'_{\ess},\lambda'_{\ess})= (\cN_{\ess},\lambda_{\ess})$.
 \end{lemma}

Let $Y, Y',Y''$ be i.i.d. unit rate Poisson processes. Then, for $s,t\ge0$, $Y(t+s)$ and $Y'(t)+Y''(s)$ have the same distribution. Thus, we have
\begin{align}\label{def_dynstoeq}
X(t)=X(0)+\sum_{y\to y'\in \phi(\cR)}Y_{y\to y'}\Big(\int_0^t\lambda_{\ess,y\to y'}(X(s))ds\Big)(y'-y).
\end{align} 
Every (weak) solution to \eqref{def_dynsto} is also a (weak) solution to \eqref{def_dynstoeq}, and vice versa.  
Consequently, the dynamics of any SRN is determined by its essential SRN as summarised in the following.

\begin{proposition}\label{lmm_seq}
Let $(\cN,\lambda)$ and $(\cN',\lambda')$ be such that $(\cN_{\ess},\lambda_{\ess})=(\cN_{\ess}',\lambda_{\ess}'')$. Then, the dynamics of $(\cN,\lambda)$ and $(\cN',\lambda')$ are equivalent, in the sense that every weak solution to  \eqref{def_dynsto} under $(\cN,\lambda)$ is also a weak solution to \eqref{def_dynsto} under $(\cN',\lambda')$, and vice versa.
\end{proposition}

\subsection{Stationary distributions}

Let $x,x'\in \Z_{\geq 0}^n$ be two states. Then,  $x$ {\bf leads to} $x'$ in $(\cN,\lambda)$, written $x\to_{\cN}x'$, if there exist reactions $y_1\ce{->}y_1',\dots, y_m\ce{->}y_m'\in \cR$, such that
\begin{enumerate}[(i)]
\item $x\geq \phi(y_1),\quad x-\phi(y_1)+\phi(y_1')\geq \phi(y_2),\quad\dots,\quad x+\sum_{i=1}^{m-1}(\phi(y_i')-\phi(y_i))\geq \phi(y_m)$.

\item $x+\sum_{i=1}^{m}(\phi(y_i')-\phi(y_i))=x'$,
\end{enumerate}
that is, the firing of  the reactions in succession will take the chain from the state $x$ to $x'$.
Condition \ref{def_cptsto} ensures that if $x\to_{\cN}x'$, then there is positive probability to jump  from $x$ to $x'$, and vice versa.

The proof of the next statement is elementary and thus omitted.
\begin{lemma}
Let $(\cN,\lambda)$ and $(\cN',\lambda')$ be SRNs with  $\cN_{\ess}=\cN'_{\ess}$.
For  $x,x'\in \Z_{\geq 0}^n$, then $x\to_{\cN} x'$ if and only if $x\to_{\cN'}x'$. As a consequence, a subset $\Gamma\subseteq \Z_{\geq 0}^n$ is an irreducible component of $(\cN,\lambda)$ if and only if it is  an irreducible component of $(\cN',\lambda')$.
\end{lemma}

The following is a consequence of Proposition \ref{lmm_seq}.
\begin{corollary}
Let $(\cN,\lambda)$ and $(\cN',\lambda')$ be SRNs such that $(\cN_{\ess},\lambda_{\ess})=(\cN'_{\ess},\lambda'_{\ess})$. Then, a probability distribution $\pi$ is a stationary distribution on an irreducible component $\Gamma$ of $(\cN,\lambda)$, if and only if $\pi$ is also a stationary distribution on $\Gamma$ of $(\cN',\lambda')$.
\end{corollary}

For complex-balanced distributions, the one-directional implication follows from \eqref{equ_kcleaved SRN}.
\begin{corollary}\label{lmm_odcb}
Let $(\cN,\lambda)$ and $(\cN',\lambda')$ be SRNs such that $(\cN',\lambda')\succeq (\cN,\lambda)$. If a probability distribution $\pi$ is a complex-balanced distribution on an irreducible component $\Gamma$ of $(\cN',\lambda')$, then $\pi$ is also a complex-balanced distribution on $\Gamma$ of $(\cN,\lambda)$.
\end{corollary}

Let $(\cN',\lambda')\succeq (\cN,\lambda)$ with projection $\psi$. For any cycle $\gamma\subseteq \cR'$, we say $\gamma$ is {\bf simple} when projected onto the digraph of $\cN$, if $\psi(\gamma)$ is a cycle of   $\cR=\psi(\cR')$. Moreover, two cycles $\gamma, \gamma'\subseteq \cR'$ are called {\bf similar} if $\psi(\gamma)=\psi(\gamma')$, when projected onto $\cN$.

\begin{theorem}\label{thm_cleaved_RN}
Let $(\cN,\lambda)$ be an SRN with a weakly reversible digraph. Then, there exists a cleaved SRN $(\cN_{\cyc},\lambda_{\cyc})$  of $(\cN,\lambda)$ with projection $\psi_{\cyc}$, such that the digraph of $\cN_{\cyc}$ consists of  
pairwise non-similar simple cycles when projected onto $\cN$, satisfying, 
\begin{enumerate}[(i)]
\setlength\itemsep{0.5em}
\item For any cycle $ \gamma \subseteq \cR_{\cyc}$, $\psi_{\cyc}( \gamma )$ is a cycle in $\cR$.
\item For any cycle $\gamma\subseteq \cR$, there exists a unique cycle $ \gamma'\subseteq \cR_{\cyc}$ such that $ \psi_{\cyc}( \gamma')=\gamma$.
\item A probability distribution $\pi$ is a complex-balanced distribution of $(\cN,\lambda)$ on some irreducible component $\Gamma$, if and of if it is one of $(\cN_{\cyc},\lambda_{\cyc})$ on $\Gamma$.
\end{enumerate}
\end{theorem}

\section{Cleaving SRNs with weakly reversible digraphs }\label{sec_cleave}

In this section, we develop an iterative procedure to show that there exists a dynamically equivalent cleaved SRN consisting of all cycles appearing in the original SRN, while preserving the complex-balanced property. This cleaving procedure enlarges the applicability of Theorem \ref{prop_sdcrt} and is key to the proof of Proposition \ref{prop_cbmcr}.

\subsection{One-node cleaving}\label{ssec_ssn}

Let $\cN=(\cC,\cR,\cS,\phi)$ be a weakly reversible RN with stochastic kinetics $\lambda$. Choose a complex $z\in \cC$ with $p_z>1$ incoming reactions.  We provide a method to construct a cleaved SRN $(\cN_1,\lambda_1)$ of $(\cN,\lambda)$ such that the complex-balanced property  of $(\cN_1,\lambda_1)$ is the same as that of $(\cN,\lambda)$, and such that $z$ is replaced by $p_z$ complexes with only one incoming reaction. Proofs are given in Section \ref{sec_prfs}.

The one-node cleaving involves two steps. In the first step, we give a precise definition of $\cN_1=(\cC_1,\cR_1,\cS,\phi_1)$ and the projection $\psi_1$, while in the second step, a kinetics is assigned to $\cN_1$. Step 1 is illustrated in Figure \ref{fig:fig1}.

{\it Step 1.} 
Order the incoming reactions of $z$ by $y_1\ce{->}z,\ \dots, \ y_{p_z}\ce{->}z$. 
 Define
\begin{align*}
\cC_1=\{y\colon y\in \cC\}\setminus\{z\}\cup \{(z,i)\colon 1\leq i\leq p_z\},
\end{align*}
and $\cR_1=\cR_1^0\cup \cR_1^{in}\cup \cR_1^{out},$  where 
\begin{align*}
\cR_1^0&=\{y\ce{->}y'\in \cR\colon y,y'\in \cC\setminus \{z\}\},\\
\cR_1^{in}&=\{y_i\ce{->}(z,i)\colon 1\leq i\leq p_z\},
\end{align*}
and $\cR_1^{out}$ is the collection of all directed edges $(z,i)\ce{->}y$ for some $i\in \{1,\dots, p_z\}$ such that there exists a cycle $\gamma$ in $\cR$ and  $y\in \cC\setminus\{z\}$  with $\{y_i\ce{->}z\ce{->}y\}\subseteq \gamma$.  By weak reversibility of $\cN$, there is at least one $i$ such that $\{y_i\ce{->}z\ce{->}y\}$ is contained in a cycle of $\cN$.

\begin{figure}
 \begin{tikzcd}[/tikz/>=stealth,/tikz/column 1/.style={column sep=2em}]
 \ch{$y_1$}\arrow[dr, ->] & & \ch{$y_1'$}\arrow[dashed,ll,->, blue, bend right=18] &  & \ch{$y_1$}\arrow[->,r] & \ch{$(z,1)$}\arrow[->,r]\arrow[->,dr] & \ch{$y_1'$}\arrow[dashed,ll,->, blue, bend right=22]\\
\ch{$\cN$:}\qquad\qquad & \ch{$z$}\arrow[ur,->] \arrow[->,r] \arrow[dr,->] & \ch{$y_2'$} \arrow[dashed,ull, ->, blue, controls={+(1.5,2.5) and +(1,1.2)}] \arrow[dashed,dll, ->, red, controls={+(1.5,-2.5) and +(1,-1.2)}] &  
&\ch{$\cN_1$:}\qquad\qquad  & & \ch{$y_2'$}\arrow[dashed,ull, ->, blue, controls={+(1.5,2.8) and +(1,1.3)}]  \arrow[dashed,dll, ->, red, controls={+(1.5,-2.8) and +(1,-1.3)}] \\
 \ch{$y_2$} \arrow[ur,->] & & \ch{$y_3'$} \arrow[dashed,ll,->, red, bend left=18] & & \ch{$y_2$}\arrow[->,r] & \ch{$(z,2)$} \arrow[->,r] \arrow[->,ur] & \ch{$y_3'$} \arrow[dashed,ll,->, red, bend left=22] &
\end{tikzcd}
\caption{ 
One-node cleaving. The complex $z$ is cleaved. A dashed edge, e.g., $y_1'$ to $y_1$, means that there exists a path directed from the initial to the terminal complex without passing through $z$. Hence, since there is a cycle containing $y_1\ce{->}z\ce{->}y_1'$ and $y_1\ce{->}z\ce{->}y_2'$, respectively, in  $\cN$, it follows that $(z,1)\ce{->}y_1'$ and $(z,1)\ce{->}y_2'$, respectively, in $\cN_1$. For the same reason, $(z,2)\ce{->}y_2'$ and $(z,1)\ce{->}y_3'$ in $\cN_1$. Primed and unprimed complexes could be the same, for example, $y_2=y_2'$. } \label{fig:fig1}
 \end{figure}

We remark  that $\{y_i\ce{->}z\ce{->}y\}\subseteq \cR$ does not imply $\{y_i\ce{->}(z,i)\ce{->}y\}\subseteq \cR_1$. For example, let $\cR=\{y_i\ce{<=>}z\ce{<=>}y\}$. Then, $\cR$ is weakly reversible and there is a closed walk $y_i\ce{->}z\ce{->}y\ce{->}z\ce{->}y_i$, including $\{y_i\ce{->}z\ce{->}y\}$. But $\{y_i\ce{->}z\ce{->}y\}$ is not  in any cycle of $\cR$, and thus $\{y_i\ce{->}(z,i)\ce{->}y\}\not\subseteq \cR_1$.

Finally, we define the labelling $\phi_1=\phi\circ \psi_1$ with $\psi_1$ the canonical projection on $\cC_1$ given by
\begin{align}\label{def_vphi1}
\psi_1(y)=\begin{cases}
y, & \text{for}\quad y\in \cC\setminus \{z\},\\
z, & \text{for}\quad y=(z,i),\quad i=1,\dots, p_z.
\end{cases}
\end{align}
 
\begin{lemma}\label{lmm_sg}
Let $\cN$ be weakly reversible, and let $\cN_1$ and $\psi_1$ be a one-node cleaved RN of $\cN$. Then, $\cN_1\succeq\cN$ with projection $\psi_1$, and $\cN_1$ is weakly reversible as well.
\end{lemma}

{\it Step 2.} 
We assign a kinetics $\lambda_1$ to $\cN_1$, such that $(\cN_1,\lambda_1)\succeq(\cN,\lambda)$ and the complex-balanced property is maintained. To complete the task, we introduce some notation. Let 
$z_1,z_2,z_3\in \cC$ be any, possibly repeated, complexes such that $\{z_1\ce{->}z_2\ce{->}z_3\}\subseteq \cR$. Denote by $\Gamma_{z_1\to z_2\to z_3}(k)$, $k\in \Z_{>0}$, the collection of closed walks in $\cR$ of the form
\begin{align*}
\gamma=\{z_1\ce{->}z_2\ce{->}z_3\ce{->}y^{(1)}\ce{->}\cdots \ce{->}y^{(k)}\ce{->}z_1\}\subseteq \cR,
\end{align*}
satisfying $\{y^{(1)},\dots, y^{(k)}\}\cap \{z_2\}=0$. For $z_1\neq z_3$, define
\[
\Gamma_{z_1\to z_2\to z_3}(0)=
\begin{cases}\{z_1\ce{->}z_2\ce{->}z_3\ce{->}z_1\}, &z_3\ce{->}z_1\in \cR\\
0, &z_3\ce{->}z_1\notin \cR,
\end{cases}
\]
and  $\Gamma_{z_1\to z_2\to z_1}(0)\coloneqq \{z_1\ce{<=>}z_2\}$. By convention,  $\Gamma_{z_1\to z_2\to z_3}(k)=0$ for  $k\in \Z_{\geq 0}$ if $\{z_1\ce{->}z_2,z_2\ce{->}z_3\}\not\subseteq \cR$.
 Define
\begin{align*}
\Gamma_{z_1\to z_2\to z_3}=\cup_{k=0}^{\infty}\Gamma_{z_1\to z_2\to z_3}(k).
\end{align*}
Furthermore, define $\rho_{z_3,z_1\to z_2}\colon\Z_{\geq 0}^n\to \R_{\geq 0}$ for all $z_1,z_2,z_3\in\cC$ and $x\in \Z_{\geq 0}^n$ by
\begin{align}\label{def_rhosto}
\rho_{z_3,z_1\to z_2}(x)=\begin{cases}\displaystyle
\frac{\lambda_{z_1\to z_2}(x+\phi(z_1)-\phi(z_3))}{\sum_{y''\colon z_1\to y''\in \cR}\lambda_{z_1\to y''}(x+\phi(z_1)-\phi(z_3))}, & z_1\ce{->}z_2\in \cR,\\
0, & z_1\ce{->}z_2\notin\cR,
\end{cases}
\end{align}
where by convention $\frac{0}{0}=0$.  
Using Condition \ref{def_cptsto}, for $z_1\ce{->}z_2\in \cR$, it holds that
$\lambda_{z_1\to z_2}(x+\phi(z_1)-\phi(z_3))>0$
 $\iff$ $x+\phi(z_1)-\phi(z_3)\geq \phi(z_1)$  
$\iff$ $x\geq \phi(z_3)$.  Thus, 
\begin{align}\label{rho_cmpt}
\rho_{z_3,z_1\to z_2}(x)>0,\quad \text{if and only if}\quad x\geq \phi(z_3).
\end{align}

\noindent Define the kinetics $\lambda_{1}$ as:
\begin{align}\label{def_lmd1out}
\lambda_{1,r}(x)=\begin{cases}
\displaystyle\lambda_{y\to y'}(x), & r=y\ce{->}y'\in\cR_1^0,\\
\displaystyle\lambda_{y_i\to z}(x), & r=y_i\ce{->} (z,i)\in \cR_1^{in},\\
\displaystyle\sum_{\gamma \in \Gamma_{y_i\to z\to y'}}\prod_{r'\in \gamma\setminus \{z\to y'\}}\rho_{z,r'}(x)\lambda_{z\to y'}(x), & r=(z,i)\ce{->} y'\in\cR_1^{out},
\end{cases}
\end{align}
for all $x\in \Z_{\geq 0}^n$. \eqref{rho_cmpt} implies that Condition \ref{def_cptsto} holds for $(\cN_1,\lambda_{1})$ as well.

\begin{lemma}\label{lmm_dswr}
Let $(\cN,\lambda)$ be a weakly reversible SRN, and $(\cN_1,\lambda_1)$ and $\psi_1$ be a one-node cleaved SRN (defined above) of $(\cN,\lambda)$. Then, $(\cN_1,\lambda_1)\succeq(\cN,\lambda)$ with projection $\psi_1$. 
\end{lemma}

By Lemma \ref{lmm_drsplt} and Proposition \ref{lmm_seq}, a stationary distribution of $(\cN,\lambda)$ is also a stationary distribution of $(\cN_1,\lambda_1)$ and vice versa. Furthermore, we have preservation of complex-balanced distributions in both directions along one-node cleavings:

\begin{lemma}\label{lmm_scb} Let $(\cN,\lambda)$ be a weakly reversible SRN, and $(\cN_1,\lambda_1)$ and $\psi_1$ be a one-node cleaved SRN (defined above) of $(\cN,\lambda)$.
Then a probability distribution $\pi$ on an irreducible component $\Gamma$ is a complex-balanced distribution of $(\cN, \lambda)$ if and only if it is a complex-balanced distribution of $(\cN_1, \lambda_1)$ on $\Gamma$.
\end{lemma}

\begin{example}
We illustrate the one-node cleaving procedure on Example \ref{exp_trag}. 
The node $A$ is cleaved. There are two incoming reactions of $A$ in $\cN$ leading to two new nodes $(A,1)$ and $(A,2)$, $\cR_1^0=\{B\ce{<=>}C\}$ and $\cR_1^{in}=\{B\ce{->}(A,1), C\ce{->}(A,2)\}$ in $\cN_1$. Since there are two cycles in $\cN$ including $B\ce{->}A$, namely $\{B\ce{->}A\ce{->}B\}$ and $\{B\ce{->}A\ce{->}C\ce{->}B\}$, then $\{(A,1)\ce{->}B,(A,1)\ce{->}C\}\subseteq \cR_0^{out}$. Similarly, we find $\{(A,2)\ce{->}B,(A,2)\ce{->}C\}\subseteq \cR_0^{out}$, and thus
$\cR_1^{out}=\{(A,1)\ce{->}B,(A,1)\ce{->}C,(A,2)\ce{->}B,(A,2)\ce{->}C\}.$ 
Consequently, the digraph of $\cN_1$ is as shown below.
\begin{equation*}
\begin{aligned}[c]
\cN:
\end{aligned}\qquad
\begin{aligned}[c]
\schemestart 
  $C$
  \arrow{<=>}[150]
  $A$
  \arrow{<=>}[210]
  $B$
  \arrow(@c3--@c1){<=>}
\schemestop
\end{aligned}
\qquad\Longrightarrow\qquad
\begin{aligned}[c]
\cN_1:
\end{aligned}\qquad
\begin{aligned}[c]
\schemestart 
  $(A,1)$
  \arrow{->}
  $C$
  \arrow{<=>}[270]
  $(A,2)$.
  \arrow{->}[180]
  $B$
  \arrow{<=>}[90]
  \arrow(@c2--@c4){<=>}
\schemestop
\end{aligned}
\end{equation*}

Concerning the kinetics of $\cN_1$, let $x=(x_A,x_B,x_C) \in \Z_{\geq 0}^3$ denote the molecular counts of the species $A$, $B$ and $C$, respectively. Using \eqref{def_lmd1out}, it suffices to calculate $\lambda_{1,(A,i)\to B}$ and $\lambda_{1,(A,i)\to C}$, $i=1,2$. Consider $(A,1)\ce{->}B\in \cR_1^{out}$. The closed walks of $\Gamma_{B\to A\to B}$ in $\cN$ are of the form
\begin{equation*}
\theta_k=\{B\ce{->}A\ce{->}B\ce{->}C\ce{->}B\ce{->}\cdots\ce{->} C\ce{->}B\},
\end{equation*}
where $k\ge0$ denotes the number of occurrences of $C\ce{->}B$. As both $B$ and $C$ have each two outgoing reactions in $\cN$, $\rho_{A,B\to C}(x)<1$ and $\rho_{A,C\to B}(x)<1$, and so
\begin{align*}
\lambda_{1,(A,1)\to B}(x)&=\sum_{k=0}^{\infty}\rho_{A,B\to A}(x)\big(\rho_{A,B\to C}(x)\rho_{A,C\to B}(x)\big)^k\lambda_{A\to B}(x) \\&=\frac{\lambda_{A\to B}(x)\rho_{A,B\to A}(x)}{1-\rho_{A,B\to C}(x)\rho_{A,C\to B}(x)}.
\end{align*}
Similarly, 
\begin{align*}
\lambda_{1,(A,1)\to C}(x)&=\frac{\lambda_{A\to C}(x)\rho_{A,B\to A}(x)\rho_{A,C\to B}(x)}{1-\rho_{A,B\to C}(x)\rho_{A,C\to B}(x)},\\
\lambda_{1,(A,2)\to B}(x)&=\frac{\lambda_{A\to B}(x)\rho_{A,B\to C}(x)\rho_{A,C\to A}(x)}{1-\rho_{A,B\to C}(x)\rho_{A,C\to B}(x)},\\
\lambda_{1,(A,2)\to C}(x)&=\frac{\lambda_{A\to C}(x)\rho_{A,C\to A}(x)}{1-\rho_{A,B\to C}(x)\rho_{A,C\to B}(x)}.
\end{align*}

By \eqref{def_rhosto}, if $x_A\geq 1$, then 
$\rho_{A,B\to A}(x)+\rho_{A,B\to C}(x)=\rho_{A,C\to A}(x)+\rho_{A,C\to B}(x)=1.$
 This implies that $\lambda_{1,(A,1)\to B}(x)+\lambda_{1,(A,2)\to B}(x)=\lambda_{A\to B}(x)$, and $\lambda_{1,(A,1)\to C}(x)+\lambda_{1,(A,2)\to C}(x)=\lambda_{A\to C}(x).$
Therefore, \eqref{equ_kcleaved SRN} holds for all $y\ce{->} y'\in \cR$ and $x\in \Z_{\geq 0}$. As a result, $(\cN_1,\lambda_1)\succeq (\cN,\lambda)$.

\end{example}

\subsection{Iteration}
 
We apply the one-node cleaving procedure iteratively until every complex has at most one incoming reaction, and the cleaved RN consists of only cycles (Lemma \ref{lmm_cycle}). 
 However, as illustrated in Figure \ref{fig:fig1}, when cleaving a complex (here, $z$), the number of incoming reactions of other complexes (here, $y_2'$) might increase. Thus, we should not expect that there is an iterative procedure based on one-node cleaving, such that the number of complexes with multiple incoming reactions is strictly decreasing.

Let $\cN=(\cC,\cR,\cS,\phi)$ be a weakly reversible RN and let $\cC'\subseteq \cC$ be the collection of complexes in $\cC$ with a single incoming reaction. Suppose that $\cC'\neq \cC$ (otherwise the RN consists of disjoint cycles, and we are done) and let $\cC''=\cC\setminus \cC'$. Write $\cN_0=(\cC_0,\cR_0,\cS,\phi_0)$ for the cleaved RN of $\cN$ with projection $\psi_0$ obtained by one-node cleaving of an arbitrary node $z\in \cC''$.  Define 
 \[
\cC_0'=\{y\in \cC_0\colon \psi_0 (y)\in \cC'\cup\{z\}\}=\cC'\cup\{y\in\cC_0 \colon \psi_0(y)=z\}\quad \mathrm{and}\quad \cC_0''= \cC_0 \setminus \cC_0'.
\]
With Figure \ref{fig:fig1} as an example, we have $\{(z,1),(z,2),y_1',y_2',y_3'\}\subseteq \cC_0'$, and $y_2'$ has two incoming reactions. Moreover, since $z\in\cC''$ and $\cC_0''=\cC''\setminus \{z\}$, then $\cC_0''$ has exactly one complex less than $\cC''$. $\cR_0$ are the reactions of the cleaved RN of $\cN$ (defined in step 1).

We next define a sequence of cleaved RNs, see Figure \ref{fig:fig2}. For $m\ge 1$,   let $\cN_m=(\cC_m,\cR_m,\cS,\phi_m)$ with projection $\psi_m$ be an RN obtained by cleaving  an element of $\cC'\cap \cC_{m-1}$ in $\cN_{m-1}$ with multiple incoming reactions (again $\cR_m$ are the reactions of the cleaved RN of $\cN_{m-1}$ as defined in step 1). Concretely, let $\psi_0^m=\psi_0\circ \dots\circ \psi_m$, 
\[
\cC_m'=\{y\in \cC_m\colon \psi_0^m(y)\in \cC'\cup\{z\}\}\subseteq\cC_m,\quad \mathrm{and}\quad \cC_m''=\cC_m\setminus \cC_m'.
\]
If all $y\in \cC'\cap\cC_{m-1}\subseteq \cC_{m-1}'$ have only one incoming reaction in $\cR_{m-1}$, then  $\cN_{m}=\cN_{m-1}$ (and $\psi_{m}=\id_{\cC_{m-1}}$). 
Hence, $\cN_{m}\succeq\cN_{m-1}$ with projection $\psi_{m}$, and $\cN_m\succeq\cN$ with projection $\psi^m_0$.  The procedure ends after $M=|\cC'|$ iterations.

\begin{figure}
\begin{floatrow}
\begin{subfigure}{0.5\textwidth}
\scalebox{0.85}{\begin{tikzpicture}
\draw node at (2.3,2) {$\cN$:};
\draw node at (4,2) {$z$};
\draw node at (6,2) {$y_1$};
\draw node at (8,2) {$y_2$};
\draw node at (5,3.7){$z'$};
\path[->] (4.3,2) edge (5.7,2);
 \path[->] (6.3,2) edge (7.7,2);
 \draw[->,blue] plot [smooth,tension=1] coordinates {(7.85,1.75) (6, 1) (4.1, 1.8)}; 
 \draw[->,dashed,red] plot [smooth,tension=1] coordinates {(7.85,2.3) (6.9, 3.7) (5.25, 3.8) };
 \path[->,red] (4.82,3.45) edge (4.05,2.2);
 \draw[thick] plot [smooth cycle] coordinates {(3.6,0) (8,-0.3) (9.8,0.2) (10,3.2) (5,4.2) (3, 2.8) };
\draw[green,thick] plot [smooth cycle] coordinates {(4.8,0.3) (8,0.4)  (10,1.5) (9,3.2) (7,3) (5.2,3.2) };
\draw node at (9.3, 1.4) {\large\color{green}$\cC'$};
\draw node at (9,0.2) {\large$\cC$};
\end{tikzpicture}}
\end{subfigure}
\begin{subfigure}{0.5\textwidth}
\scalebox{0.85}{\begin{tikzpicture}
\draw node at (2.3,2) {$\cN_0$:};
\draw node at (4,2.5) {$(z,2)$};
\draw node at (6.2,2) {$y_1$};
\draw node at (4,1.5) {$(z,1)$};
\draw node at (8,2) {$y_2$};
\draw node at (5,3.7){$z'$};
\path[->] (4.55,1.5) edge (5.9,1.8);
 \path[->] (6.5,2) edge (7.7,2);
 \draw[->,blue] plot [smooth,tension=1] coordinates {(7.85,1.7) (6, 0.5) (4.2, 1.2)};
 \path[->] (4.55,2.5) edge (5.9,2.1);
 \draw[->,dashed,red] plot [smooth,tension=1] coordinates {(7.85,2.3) (6.9, 3.7) (5.27, 3.8) };
 \path[->,red] (4.75,3.5) edge (4.1,2.85);
 \draw[thick] plot [smooth cycle] coordinates {(3.6,0) (8,-0.3) (9.8,0.2) (10,3.2) (5,4.2) (3, 2.8) };
\draw[green,thick] plot [smooth cycle] coordinates {(4,0.3) (8,0.4)  (10,1.5) (9,3.2) (7,3) (3.5,3.2) };
\draw node at (9.3, 1.4) {\large\color{green}$\cC_0'$};
\draw node at (9,0.2) {\large$\cC_0$};
\end{tikzpicture}}
\end{subfigure}
\end{floatrow}
\begin{floatrow}
\begin{subfigure}{0.5\textwidth}
\scalebox{0.85}{\begin{tikzpicture}
\draw node at (2.3,2) {$\cN_1$:};
\draw node at (4,2.5) {$(z,2)$};
\draw node at (6.2,2.5) {$(y_1,2)$};
\draw node at (4,1.5) {$(z,1)$};
\draw node at (6.2,1.5) {$(y_1,1)$};
\draw node at (8.2,2) {$y_3$};
\draw node at (5,3.7){$z'$};
\path[->] (4.6,1.5) edge (5.55,1.5);
 \path[->] (6.85,1.55) edge (7.9,1.9);
 \draw[->,blue] plot [smooth,tension=1] coordinates {(8.05,1.7) (6, 0.5) (4.1, 1.25)};
\path[->] (4.6,2.5) edge (5.55,2.5);
 \path[->] (6.85,2.45) edge (7.9,2.1);
 \draw[->,dashed,red] plot [smooth,tension=1] coordinates {(8.05,2.3) (6.9, 3.7) (5.2, 3.8) };
 \path[->,red] (4.8,3.5) edge (4,2.8);
 \draw[thick] plot [smooth cycle] coordinates {(3.6,0) (8,-0.3) (9.8,0.2) (10,3.2) (5,4.2) (3, 2.8) };
\draw[green,thick] plot [smooth cycle] coordinates {(4,0.3) (8,0.4)  (10,1.5) (9,3.2) (7,3) (3.5,3.2) };
\draw node at (9.3, 1.4) {\large\color{green}$\cC_2'$};
\draw node at (9,0.2) {\large$\cC_2$};
\end{tikzpicture}}
\end{subfigure}
\begin{subfigure}{0.5\textwidth}
\scalebox{0.85}{\begin{tikzpicture}
\draw node at (2.3,2) {$\cN_2$:};
\draw node at (4,2.5) {$(z,2)$};
\draw node at (6.4,2.5) {$(y_1,2)$};
\draw node at (8.8,2.5) {$(y_3,2)$};
\draw node at (4,1.8) {$(z,1)$};
\draw node at (6.4,1.8) {$(y_1,1)$};
\draw node at (5.2,0.5) {$(y_3,1)$};
\draw node at (5,3.7){$z'$};
\path[->] (4.6,1.8) edge (5.7,1.8);
 \path[->] (6.3,1.5) edge (5.3,0.8);
 \path[->,blue] (4.1,1.5) edge (5.1,0.8);
\path[->] (4.6,2.5) edge (5.7,2.5);
 \path[->] (7.1,2.5) edge (8.1,2.5);
 \draw[->,dashed,red] plot [smooth,tension=1] coordinates {(8.7,2.85) (7.2, 3.8) (5.3, 3.75) };  
 \path[->,red] (4.78,3.48) edge (4,2.8);
 \draw[thick] plot [smooth cycle] coordinates {(3.6,0) (8,-0.3) (9.8,0.2) (10,3.2) (5,4.2) (3, 2.8) };
\draw[green,thick] plot [smooth cycle] coordinates {(4,0.3) (8,0.4)  (10,1.5) (9,3.2) (7,3) (3.5,3.2) };
\draw node at (9.3, 1.4) {\large\color{green}$\cC'_3$};
\draw node at (9,0.2) {\large$\cC_3$};
\end{tikzpicture}}
\end{subfigure}
\end{floatrow}\caption{}\label{fig:fig2}
\end{figure}

\begin{lemma}\label{lmm_cfnt}
Every complex in $\cC'_M\subseteq\cC_M$ has only one incoming reaction in $\cR_M$.
\end{lemma}

After completing the $M$-th  iteration, we obtain a cleaved RN $\cN_M=(\cC_M,\cR_M,\cS,\phi_M)$ of $\cN$ with projection $\psi_M$, such that each complex in $\cC_M'\subseteq\cC_M$ has only one incoming reaction, and $\cC_M'' $ has one fewer complexes than $ \cC''$, namely $ \cC_M''=\cC''\setminus\{z\}$. However, the number of incoming reactions of a complex $y\in \cC_M''$ might be different from the corresponding number of incoming reactions of the complex $\psi^M_0(y)=y\in \cC''$ in $\cR$.

By repeating this procedure for another complex $z'\in \cC_M''$ and so forth, we eventually obtain, after finitely many iterations, a cleaved SRN $(\cN_{\cyc},\lambda_{\cyc})$ with projection $\psi_{\cyc}$ on $\cN$. Every complex in the cleaved SRN has only one incoming reaction. Hence, the cleaved SRN consists of disjoint cycles (Lemma \ref{lmm_cycle}). Furthermore, it has the complex-balanced property if and only if $(\cN,\lambda)$ fulfils it (Lemma \ref{lmm_scb}).

\subsection{Completion}\label{sec_cplt}

We modify the cleaved SRN $(\cN_{\cyc},\lambda_{\cyc})$ to obtain another cleaved SRN of $(\cN,\lambda)$ without non-simple cycles and similar cycles when projected onto $\cN$. The modification includes two steps. In the first step, we cut and adhere non-simple cycles, and in the second step, we combine similar cycles (for definitions see just before Theorem \ref{thm_cleaved_RN}).

Suppose there exists a cycle $\gamma\subseteq \cR_{\cyc}$ that is not simple when projected onto $\cN$. Then, it is of the form
\begin{align*}
\gamma=\big\{&y_0\ce{->}y_1\ce{->}\cdots \ce{->}y_k\ce{->}y_0'\ce{->}y_{k+1}\ce{->}\cdots\ce{->} y_{k+k'}\ce{->}y_0\big\},
\end{align*}
where $y_0\neq y_0'$ and $\psi_{\cyc}(y_0)=\psi_{\cyc}(y_0')$. 
We cut this cycle at $y_0$ and $y_0'$, then adhere each piece with its end node. Thus, we get two cycles,
\begin{align*}
\gamma_1&=\{y_0\ce{->}y_1\ce{->}\cdots \ce{->}y_k\ce{->}y_0\},\\
\gamma_2&=\{y_0'\ce{->}y_{k+1}\ce{->}\cdots \ce{->}y_{k+k'}\ce{->}y_0'\}.
\end{align*}
In this way, we obtain a new cleaved RN $\cN'_{\cyc}=(\cC'_{\cyc},\cR'_{\cyc},\cS,\phi'_{\cyc})$ of $\cN$ with projection $\psi'_{\cyc}=\psi_{\cyc}$, where $\cC'_{\cyc}=\cC_{\cyc}$,
\begin{align*}
\cR'_{\cyc}=&\big(\cR_{\cyc}\setminus \{y_k\ce{->}y_0',y_{k+k'}\ce{->}y_0\}\big)\cup \{y_k\ce{->}y_0,y_{k+k'}\ce{->}y_0'\}.
\end{align*}
It is natural to assign a kinetics $\lambda'_{\cyc}$ to $\cN'_{\cyc}$ by keeping  the same kinetics for the reactions also appearing in $\cR_{\cyc}$, and letting
\[
\lambda'_{\cyc,y_k\to y_0}=\lambda_{\cyc,y_k\to y_0'},\quad \lambda'_{\cyc,y_{k+k'}\to y_0'}=\lambda_{\cyc,y_{k+k'}\to y_0}.
\] 
Then,  $(\cN'_{\cyc},\lambda'_{\cyc})\succeq(\cN,\lambda)$ with projection $\psi_{\cyc}$, such that the complex-balancedness remains. Note that $(\cN_{\cyc},\lambda_{\cyc})$ and $(\cN'_{\cyc},\lambda'_{\cyc})$ may not be related by $\succeq$.

The `cut-adhere' process can be accomplished in finitely many steps until every cycle is simple when projected onto $\cN$.  By abuse of notation,
the final cleaved SNR is also denoted by $(\cN_{\cyc},\lambda_{\cyc})$ with projection $\psi_{\cyc}$. 

In the second step, we combine similar cycles. Suppose there are two similar cycles $\gamma_1,\gamma_2\subseteq\cR_{\cyc}$ when projected onto $\cN$, that is, $\psi_{\cyc}(\gamma_1)=\psi_{\cyc}(\gamma_2).$ We simply remove $\gamma_2$ and sum the kinetics of each reaction in $\gamma_2$ to the corresponding  reaction in $\gamma_1$. More precisely, suppose 
\begin{align*}
\gamma_1&=\{y_1\ce{->}\cdots \ce{->} y_k\ce{->}y_1\},\\
\gamma_2&=\{y_1'\ce{->}\cdots \ce{->} y_k'\ce{->}y_1'\},
\end{align*}
with $y_i\neq y_i'$ and $\psi_{\cyc}(y_i)=\psi_{\cyc}(y_i')$ for all $i=1,\dots, k$. Then, we construct a new cleaved RN $(\cN_{\cyc}',\lambda_{\cyc}')$ of $(\cN,\lambda)$ with $\psi_{\cyc}'$ being a restriction of $\psi_{\cyc}$ on $\cC_{\cyc}'=\cC_{\cyc}\setminus \{(y_j,i_j')\colon 1\leq j\leq k\}$, where $\cR_{\cyc}'=\cR_{\cyc}\setminus \gamma_2$, the labelling $\phi_{\cyc}'$ is a restriction of $\phi_{\cyc}$ on $\cC_{\cyc}'$, and the kinetics $\lambda_{\cyc}'$ is defined as follows,
\[
\lambda_{\cyc,r}'=\begin{cases}
\lambda_{\cyc,r}, & r\in \cR_{\cyc}\setminus (\gamma_1\cup \gamma_2),\\
\lambda_{\cyc,y_j\to y_{j+1}}+\lambda_{\cyc,y_j'\to y_{j+1}'}, &r=y_j\ce{->}y_{j+1}\in \gamma_1.
\end{cases}
\]
Then $(\cN_{\cyc}',\lambda_{\cyc}')\succeq(\cN,\lambda)$ with projection $\psi_{\cyc}'$, and $(\cN_{\cyc}',\lambda_{\cyc}')$ fulfils the complex-balancedness if and only if $(\cN,\lambda)$ fulfils it. Here, $(\cN_{\cyc},\lambda_{\cyc})\succeq(\cN_{\cyc}',\lambda_{\cyc}')$.

This process can be iterated finitely many times until all disjoint cycles are non-similar when projected onto $\cN$. By abuse of notation, the resulting cleaved SRN of $(\cN,\lambda)$ is also denoted by $(\cN_{\cyc},\lambda_{\cyc})$ with projection $\psi_{\cyc}$.

\section{Proofs}\label{sec_prfs}

 \subsection{Proof of Theorem \ref{prop_sdcrt}}  
 The proof follows the idea of \cite[Theorem 4.1]{anderson2}. By definition, to show that $\pi$  
is a complex-balanced distribution, it suffices to verify that for any complex $\eta\in \cC$ and any $x\in \Gamma$, the following holds
\begin{align}\label{psd1}
g(x)\sum_{y' \colon \eta\to y'\in \cR}\lambda_{\eta\to y'}(x)=\sum_{y\colon y\to \eta\in \cR}g\big(x+\phi(y)-\phi(\eta)\big)\lambda_{y\to \eta}\big(x+\phi(y)-\phi(\eta)\big).
\end{align}
Due to Condition \ref{def_cptsto}, we only need to prove \eqref{psd1}   assuming $x\geq \phi(\eta)$. Note that for any $\eta\in \cC$, all reactions such that $\eta$ is a reactant or product are in the same connected component. 
Assume $\eta\in\cL_k$. Then,  \eqref{sdcrt} yields   for all $x\in \Gamma$ and $x\geq \phi(y)$, 
\begin{align}\label{psd2}
&\sum_{y\colon y\to \eta\in \cR }\lambda_{y\to \eta}\big(x+\phi(y)-\phi(\eta)\big) g\big(x+\phi(y)-\phi(\eta)\big)=\frac{\sum_{ y\colon y\to \eta\in \cR }\kappa_{y\to \eta}}{m_{k}(x-\phi(\eta))},
\end{align}
and
\begin{align}\label{psd3}
\sum_{ y'\colon \eta\to y'\in \cR}\lambda_{\eta\to y'}(x)g(x)=\frac{\sum_{y'\colon \eta\to y'\in \cR }\kappa_{\eta \to y'}}{m_{k}(x-\phi(\eta))}.
\end{align}
Then, equality \eqref{psd1} follows from \eqref{sdcrt0}, \eqref{psd2} and \eqref{psd3}.  \qed

\subsection{Proofs of Proposition \ref{prop_cbmcr}}

  As a consequence of Theorem \ref{prop_sdcrt}, we only need to show one direction. Suppose that $\pi$ is a complex-balanced distribution of $(\cN,\lambda)$ on   $\Gamma$.

Recall the assumption that $\lambda_{y\to y'}=\alpha_{y\to y'}\lambda_y $ on $\Gamma$ for all $y\ce{->} y'$. For any $\eta\in \cC$ and $r\in \cR$, the function $\rho_{\eta,r}$  in \eqref{def_rhosto} is a positive constant on $\{x\in \Gamma\colon x\geq \phi(\eta)\}$. Therefore, $\lambda_{1,r}$  in \eqref{def_lmd1out} fulfils $\lambda_{1,r}(x)=c(r)\lambda_{\psi_1(r)}(x)$ for some constant $c(r)$. After iteration and completion as in Section \ref{sec_cleave}, we  find a cleaved SRN $(\cN_{\cyc},\lambda_{\cyc})$ of $(\cN,\lambda)$ with projection $\psi_{\cyc}$, such that 
\begin{align*}
\lambda_{\cyc,r}(x)=c(r)\lambda_{\psi_{\cyc}(r)}(x),
\end{align*}
for all $r\in \cR_{\cyc}$ and $x\in \Gamma$ with positive constants $\{c(r'), r'\in \cR_{\cyc}\}$. Choose any $r=(y,i)\ce{->}(y',i')\in \cR_{\cyc}$, where $(y,i),(y',i')\in \cC_{\cyc}$, such that $\psi_{\cyc}(y,i)=y$ and $\psi_{\cyc}(y',i')=y'$. Suppose that $r$ is in the $k$-th  
connected component (cycle) of $\cN_{\cyc}$.
Using Theorem \ref{thm_cleaved_RN} and Proposition \ref{prop_cbmpf}, since $\pi$ is a complex-balanced distribution of $(\cN,\lambda)$ (and thus of $(\cN_{\cyc},\lambda_{\cyc})$) on $\Gamma$, we have
\begin{align}\label{for_pisp}
\pi(x)=\big[\lambda_{\cyc,r}(x)m_{\cyc,k}(x-\phi_{\cyc}(y,i))\big]^{-1}=\big[c(r)\lambda_{y\to y'}(x)m_{\cyc,k}(x-\phi(y))\big]^{-1},
\end{align}
for all $x\in \Gamma$. Then, the proposition follows if we can show that the ratio $m_{j_1,\cyc}/m_{j_2,\cyc}$ is a constant on $\Gamma_j$ (see \eqref{def_gmmk}) for any indices $j_1$, $j_2$ and $j$, such that the $j_1$-th and $j_2$-th cycles in $\cN_{\cyc}$ are both included in the $j$-th connected component
when projected onto $\cN$. The following lemma follows from weak reversibility and the proof is omitted.
 \begin{lemma}\label{lmm_gmmj}
Let $\cN$ be a weakly reversible RN consisting of connected components $\cL_1,\dots, \cL_l$. Suppose $\Gamma \subseteq \R^n$ is an irreducible component of $\cN$. For any $j\in \{1,\dots,l\}$, let $\Gamma_j$ be given as in \eqref{def_gmmk}. Then, $\Gamma_j = \{x -\phi(y)\colon x\in \Gamma\} \cap \Z_{\geq 0}^n$, where $y$ is an arbitrary complex in $\cL_k$.
 \end{lemma}

For $\iota=1,2$, let $r_{\iota}=(y_{\iota},i_{\iota})\ce{->}(y_{\iota}',i_{\iota}')$ be in the $j_{\iota}$-th cycle of $\cN_{\cyc}$, written as $r_{\iota}\in \cL_{\cyc,j_{\iota}}$. By convention, we assume $\psi_{\cyc}(r_{\iota})=y_{\iota}\ce{->}y_{\iota}'$. 
Furthermore, suppose that $\psi_{\cyc}(r_1)$ and $\psi_{\cyc}(r_2)$ are both in the $j$-th connected component 
 of $\cN$.

 {\bf Case 1) } Suppose $y_1=y_2$. By assumption, $\lambda_{y_1\to y_1'}/\lambda_{y_2\to y_2'}=\alpha_{y_1\to y_1'}/\alpha_{y_2\to y_2'}$ is a positive constant on $\Gamma$. Moreover, due to equation \eqref{for_pisp}, it holds for every $x\in \Gamma$ with $x - \phi (y_1) \in \Z_{\geq 0}^n$, 
\[
1=\frac{\pi(x)}{\pi(x)}=\frac{c(r_1)\lambda_{y_1\to y_1'}(x)m_{\cyc,j_1}(x-\phi(y_1))}{c(r_2)\lambda_{y_2\to y_2'}(x)m_{\cyc,j_2}(x-\phi(y_2))}=\frac{c(r_1)\alpha_{y_1\to y_1'}m_{\cyc,j_1}(x-\phi(y_1))}{c(r_2)\alpha_{y_2\to y_2'}m_{\cyc,j_2}(x-\phi(y_2))}.
\]
By assumption $y_1 = y_2$, and performing a change of variable $z = x-\phi(y_1) = x-\phi(y_2)$, we get
\begin{align}\label{eq_cs1}
\frac{m_{\cyc,j_1}(z)}{m_{\cyc,j_2}(z)}=\frac{c(r_2)\alpha_{y_2\to y_2'}}{c(r_1)\alpha_{y_1\to y_1'}},
\end{align}
is a positive constant, for every $z \in \Gamma_j$ such that $z=x - \phi(y_1)$ with some $x\in\Gamma$. Taking Lemma \ref{lmm_gmmj} into account, the identity \eqref{eq_cs1} holds for all $z\in \Gamma_j$.

{\bf Case 2)} Suppose that $y_1'=y_2$. Consider reaction $r_2$ and the outgoing reaction of $(y_1',i_1')$: $r_1'=(y_1',i_1')\ce{->}(y_1'',i_1'')\in\cL_{\cyc,j_1}$. Then, by application of Case 1, we immediately get that 
\begin{align*}
\frac{m_{\cyc,j_1}(z)
}{m_{\cyc,j_2}(z)
}=\frac{c(r_2)\alpha_{y_2\to y_2'}}{c(r_1')\alpha_{y_1'\to y_1''}},
\end{align*}
for all $z = x - \phi(y_2) = x - \phi(y_1')$ with $x\in \Gamma$, and thus for all $z\in\Gamma_j$.

{\bf Case 3)} Suppose $y_1=y_2'$. One can verify that $m_{\cyc,j_1}(z)
/m_{\cyc,j_2}(z)
$ is a positive constant on $\Gamma_j$ for every $z\in \Gamma_j$ 
following the same lines as in Case 2.

{\bf Case 4)} Suppose $y_1'=y_2'$. Consider the reactions $r_1'=(y_1',i_1')\ce{->}(y_1'',i_1'')\in \cL_{\cyc,j_1}$ and $r_2'=(y_2',i_2')\ce{->}(y_2'',i_2'')\in \cL_{\cyc,j_2}$, where $\psi_{\cyc}(y_{\iota}'',i_{\iota}'')=y_{\iota}''$ for $\iota=1,2$. Using Case 1, 
\[
\frac{m_{\cyc,j_1}(z)
}{m_{\cyc,j_2}(z)
}=\frac{c(r_2')\alpha_{y_2'\to y_2''}}{c(r_1')\alpha_{y_1'\to y_1''}},
\]
for all $z = x - \phi (y_1') = x - \phi(y_2')$ with $x\in \Gamma$, and thus for all $z\in \Gamma_j$.

{\bf Remaining cases.} Since $\psi_{\cyc}(r_1)$ and $\psi_{\cyc}(r_2)$ are both in the $j$-th connected component in $\cN$, we can find $y^{(1)},\dots, y^{(k)}\in \cC$ with $y^{(1)}=y_1$ and $y^{(k)}=y_2$, such that for all $i=1,\dots, k-1$, either $y^{(i)}\ce{->}y^{(i+1)}$ or $y^{(i+1)}\ce{->}y^{(i)}$ in $\cR$. This yields that for some indexes $q_1,q_1',\dots, q_k, q_k'$, we have $(y^{(i)},q_i)\ce{->}(y^{(i+1)},q_i')$ or $(y^{(i+1)},q_i)\ce{->}(y^{(i)},q_i')$ in the $j_i'$-th connected component in $\cN_{\cyc}$ for all $i=1,\dots, k-1$. It follows from Cases 1-4, that $m_{\cyc,j_1}(z)/m_{\cyc,j_1'}(z) = c_0$, $m_{\cyc,j_i'}(z)/m_{\cyc,j_{i+1}'}(z) = c_i$, $i=1,\dots, k-2$, and $m_{\cyc,j_{k-1}'}(z)/m_{\cyc,j_2}(z) = c_{k-1}$ for all $z\in \Gamma_j$ with some positive constants $c_0,\dots, c_{k-1}$. Thus, $m_{\cyc,j_1}(z)/m_{\cyc,j_2} (z) = c_k$ with some positive constant $c_k$ for all $z \in\Gamma_j$.

Therefore, $\pi$ can be written as in the form  \eqref{for_cpi} with appropriate positive constants $\kappa_{ y\to y'}$. Moreover,  \eqref{sdcrt0} is a direct result of  \eqref{def_vcbm} and \eqref{for_cpi}. The proof is complete. \qed

\subsection{Proof of Theorem \ref{prop_cbmpf}}
 Under the given conditions, every complex has exactly one incoming reaction and one outgoing reaction. By Theorem \ref{prop_sdcrt}, it is enough to show that if there is a complex-balanced stationary distribution, then it   satisfies \eqref{for_pidcp}. Without loss of generality, assume $\ell=1$. Then, there exists an integer $p\geq 2$, such that
\[
\cR=\big\{y_i\ce{->}y_{i+1}\colon i=1,\dots, p; \ y_k\neq y_j, 1\leq k<j\leq p; \ y_{p+1}=y_1\big\}.
\] 
Since $\pi$ is a complex-balanced distribution on $\Gamma$, then, by definition, we have,
\begin{align}\label{for_cbm1cl}
\pi(x)\lambda_{y_i\to y_{i+1}}(x)=\pi\big(x+ \phi(y_{i-1})- \phi(y_i)\big)\lambda_{y_{i-1}\to y_i}\big(x+ \phi(y_{i-1})- \phi(y_i)\big)
\end{align}
for all $i=1,\dots, p$ (by convention $y_0=y_{p+1}$) and $x+ \phi(y_1)- \phi(y_2)\in \Gamma$. We define the function $m$ as follows. For all $x\in \N_0^n$ such that $x+ \phi(y_1)\in \Gamma$, we let
\[
m(x)=\big[\pi(x+ \phi(y_1))\lambda_{y_1\to y_2}(x+ \phi(y_1))\big]^{-1},
\]
Thus, $\pi(x)=[\lambda_{y_1\to y_2}(x)m(x- \phi(y_1))]^{-1}$ if $x\in \Gamma$ with $x\geq y_1$. On the other hand,  \eqref{for_cbm1cl} yields
 \begin{align}\label{for_cbm1cl1}
\pi(x)\lambda_{y_2\to y_{3}}(x)=\pi\big(x+ \phi(y_{1})- \phi(y_2)\big)\lambda_{y_{1}\to y_2}\big(x+ \phi(y_{1})- \phi(y_2)\big)=m(x- \phi(y_2))^{-1},
\end{align}
for all $x\in \Gamma$ with $x+y_1-y_2\in \Gamma$ and $x\geq y_2$. Note that for all $x\in \Gamma$ with $x\geq \phi(y_2)$, we have $x- \phi(y_2)+ \phi(y_3)\in \Gamma$. Thus $x- \phi(y_2)+ \phi(y_4)\in \Gamma$ as well. By iteration and the fact that $y_{p+1}=y_1$, it follows that $x- \phi(y_2)+ \phi(y_1)\in \Gamma$. Therefore, \eqref{for_cbm1cl1} holds for $x\in \Gamma$ with $x\geq y_2$. This implies that $\pi(x)=[\lambda_{y_2\to y_{3}}(x)m(x- \phi(y_2))]^{-1}$, for all $x\in \Gamma$ with $x\geq \phi(y_2)$. Finally, by iteration, \eqref{for_pidcp} holds for all $y_i\to y_{i+1}$, $i=1,\dots, p$. The proof is complete. \qed

\subsection{Proof of Proposition \ref{prop_DB}}
   Suppose that \eqref{eq_prop_db} holds. Then, for any $x\in \Gamma$ and $y \ce{<=>} y' \in \cR$,
  \begin{align*}
    \pi(x) \lambda_{y\to y'}(x) = & m_{y \to y'} (x - \phi(y))^{-1} \\
    = & m_{y' \to y} (x - \phi(y))^{-1}
    = \pi (x + \phi (y') - \phi (y)) \lambda_{y'\to y} \big(x+\phi(y')-\phi(y) \big)
  \end{align*}
Consequently, $\pi$ is a detailed-balanced distribution for $(\cN,\lambda)$ on $\Gamma$.

Oppositely, suppose $(\cN,\lambda)$ is detailed-balanced on $\Gamma$ with distribution $\pi$. For any $y \ce{<=>} y' \in \cR$, define $m_{y \to y'}$ and $m_{y' \to y'}$ on $\Gamma_k$ by
\begin{align*}
  m_{y \to y'} (x) &\coloneqq \big[ \lambda_{y \to y'} (x + \phi (y)) \pi (x + \phi (y)) \big]^{-1},\\
  m_{y' \to y} (x) &\coloneqq \big[ \lambda_{y' \to y} (x + \phi (y')) \pi (x + \phi (y')) \big]^{-1}
\end{align*}
Then, by definition of detailed-balanced distribution, we have $m_{y \to y'} = m_{y' \to y}$ on $\Gamma_k$, and we are done. \qed

\subsection{Proof of Theorem \ref{thm_cleaved_RN}}
 Let $(\cN_{\cyc},\lambda_{\cyc})$ be the cleaved SRN obtained by the iterative one-node cleaving procedure (Sections \ref{ssec_ssn}-\ref{sec_cplt}). Then, by Lemmas \ref{lmm_cycle}, \ref{lmm_sg} and \ref{lmm_cfnt}, $\cN_{\cyc}$, we see that the digraph of $\cN_{\cyc}$ consists of disjoint cycles that are pairwise non-similar simple cycles when projected onto $\cN$. It suffices to check (i)-(iii) for $(\cN_{\cyc},\lambda_{\cyc})$.

First, (i) is a direct consequence of $\cN_{\cyc}\succeq\cN$ and the fact that every cycle in $\cN_{\cyc}$ is simple when projected onto $\cN$. Next, denote by $(\cN_{\cyc}^*,\lambda_{\cyc}^*)$ be the cleaved SRN before completion in Section \ref{sec_cplt}. Then, due to Lemma \ref{lmm_scb} and an iteration argument, (iii) holds for $(\cN_{\cyc}^*,\lambda_{\cyc}^*)$. The `cut-adhere' process does not affect the validity of (iii). We need to show that (iii) still holds after the combination process, which seems wrong as in Example \ref{cexp}. Still denote by $(\cN_{\cyc}^*,\lambda_{\cyc}^*)$ the cleaved SRN after `cut-adhere' before combination. Then, $(\cN_{\cyc}^*,\lambda_{\cyc}^*)\succeq (\cN_{\cyc},\lambda) \succeq (\cN,\lambda)$, and thus (iii) follows as a result of Corollary \ref{lmm_odcb}. 

Last, we need to prove (ii). Let 
\[
\gamma=\{y_{1}\ce{->}y_{2}\ce{->}\cdots\ce{->}y_{m}\ce{->}y_{1}\}
\]
be a cycle in $\cR$. Let $\cN_1=(\cC_1,\cR_1,\cS,\phi_1)$ be the cleaved RN of $\cN$ with projection $\psi_1$ due to one-node cleaving of some $z\in \cC$. If $z\notin \{y_1,\dots, y_m\}$, then every reaction of the cycle is also in $\cR_1$. On the other hand, without loss of generality, assume that $z=y_{1}$. Then, there exists some index $i\in \{1,\dots, p_z\}$, such that $y_{m}\ce{->}(z,i)$, and by definition of $\cR_1^{out}$, we have  $(z,i)\ce{->}y_2\in \cR_1$ as well. In other words, there exists a cycle 
\[
\gamma_1=\{(z,i)\ce{->}y_{2}\ce{->}\cdots\ce{->}y_{m}\ce{->}(z,i)\}\in \cR_1.
\]
Therefore, there exists a cycle $\gamma_1\in \cR_1$, such that $\psi_1(\gamma_1)=\gamma$ in any case. By iteration, there exists a cycle $\gamma_{\cyc}^*\subseteq \cR_{\cyc}^*$, such that $\psi_{\cyc}^*(\gamma_{\cyc}^*)=\gamma$, where $\cN_{\cyc}^*=(\cC_{\cyc}^*,\cR_{\cyc}^*,\cS,\phi_{\cyc}^*)$ denotes the cleaved RN of $\cN$ with projection $\psi_{\cyc}^*$, before the completion step in Section \ref{sec_cplt}. Since $\gamma_{\cyc}^*$ is simple, it will not be affected in the `cut-adhere' process. Finally, in the combination process of Section \ref{sec_cplt}, the cycle $\gamma_{\cyc}^*$ may be `absorbed' by other similar cycles when projected onto $\cN$. However, it does not influence the validity of property (ii). The proof is complete. \qed

\subsection{Proof of Lemma \ref{lmm_sg}}
We first prove that $\cN_1\succeq\cN$ with projection $\psi_1$. By definition, it suffices to show that $(\cC,\cR)=(\psi_1(\cC_1),\psi_1(\cR_1))$. In fact, due to \eqref{def_vphi1}, we have 
\[
\psi_1(\cC_1)=(\cC\setminus \{z\})\cup \{z\}=\cC.
\]
By definition of $\cR_1$, we have $\psi_1(\cR_1)\subseteq \cR$. To prove the reverse inclusion, we decompose $\cR=\cR^0\cup\cR^{in}\cup\cR^{out}$, where $\cR^0$ consists of reactions whose reactant and product are both in $\cC\setminus\{z\}$, and $\cR^{in}$ and $\cR^{out}$ consist of the incoming and outgoing reactions of $z$ in $\cR$, respectively. Then,  $\psi_1(\cR^0_1)=\cR_1^0=\cR^0$ and $\psi_1(\cR_1^{in})=\{y_i\ce{->}z\colon 1\leq i\leq p_z\}=\cR^{in}$. Recall that $\cN$ is weakly reversible. Thus, for every $y\in \cC$ such that $z\ce{->}y\in \cR$, there exists a cycle containing $z\ce{->}y$, and the incoming reaction of $z$ in this cycle is $y_j\ce{->}z$ for some $1\leq j\leq p_z$. Then, $(z,j)\ce{->}y\in \cR_1^{out}$, and thus $z\ce{->}y\in \psi_1( \cR_1^{out})$. This implies $\cR^{out}\subseteq \psi_1(\cR_1^{out})$. Thus, $\cN_1\succeq\cN$ with $\psi_1$.

Next, we show weak reversibility of $\cN_1$. Suppose that $y\ce{->}y'\in \cR_1^0$. Then, $y\ce{->}y'\in \cR$. By  weak reversibility of $\cR$, there exists a cycle $\gamma\subseteq \cR$ containing $y\ce{->}y'$. If $z\notin\gamma$, then $\gamma \subseteq \cR_1^0$, and we are done. Otherwise, suppose $z\in \gamma$, then there exist $i\in \{1,\dots, p_z\}$ and $y'\in \cC\setminus \{z\}$, such that $\{y_i\ce{->}z\ce{->}y'\}\subseteq \gamma\subseteq \cR$. As a consequence, $\{y_i\ce{->}(z,i)\ce{->}y'\}\subseteq \cR_1$. Replacing $z$ by $(z,i)$ in $\gamma$, we get a new cycle $\gamma' \subseteq \cR_1$. For reactions in $\cR_1^{in}$ or $\cR_1^{out}$, the same idea is applicable and the details are omitted.
The proof is complete. \qed

\subsection{Proof of Lemma \ref{lmm_dswr}}

Due to Lemma \ref{lmm_sg}, it suffices to show that $\lambda_1$ satisfies \eqref{equ_kcleaved SRN}. In fact, by definition of $\lambda_1$, we only need to prove that
\begin{align}\label{dswr1}
\lambda_{z\to y'}(x)=\sum_{i=1}^{p_z}\lambda_{1,(z,i)\to y'}(x)
\end{align}
for  $y'\in\cC$ with $z\ce{->}y'\in \cR$ and $x\in \Z_{\geq 0}^n$. Using \eqref{def_lmd1out} and Condition \ref{def_cptsto} on $(\cN,\lambda)$, then \eqref{dswr1} is equivalent to 
\begin{align}\label{dswr2}
\1_{\{x'\colon x'\geq \phi(z)\}}(x)=\sum_{i=1}^{p_z}\sum_{\gamma \in \Gamma_{y_i\to z\to y'}}\prod_{r\in \gamma\setminus \{z\to y'\}}\rho_{z,r}(x),
\end{align}
which is what we will prove. First, if $x\not\geq \phi(z)$, by \eqref{rho_cmpt}, both sides of \eqref{dswr2} are equal to zero.

 Hence assume that $x\geq \phi(z)$.  Let $\cX$ be the set of all complexes in $\cC$ that are in the same connected component 
as $z$. Then, weak reversibility, Condition \ref{def_cptsto} and \eqref{def_rhosto} imply that for any $z_1,z_2\in \cX$,
\begin{enumerate}[(i)]
\setlength\itemsep{0.5em}
\item $\phi(z_1)\to_{\cN}\phi(z_2)$.

\item $\rho_{z,z_1\to z_2}(x)>0$ if and only if $z_1\ce{->}z_2\in \cR$.

\item $\displaystyle\sum_{z'\in \cX}\rho_{z,z_1\to z'}(x)=\sum_{z'\colon z_1\to z'\in \cR} \rho_{z,z_1\to z'}(x)=1$.
\end{enumerate}
This observation allows us to define a discrete time Markov chain (DTMC)  on $\cX$ with transition probability $ \mathbb{P}_{z_1}(z_2)=\rho_{z,z_1\to z_2}(x)$
for all $z_1,z_2\in \cX$. Moreover, the chain is irreducible with finite state space. Therefore, it follows from \cite[Theorems 1.5.6 and 1.5.7]{norris98} that the chain is recurrent, and thus,
\begin{align*}
\mathbb{P}_{y'}(\tau_z<\infty)=\sum_{i=1}^{p_z}\sum_{\gamma \in \Gamma_{y_i\to z\to y'}}\prod_{r\in \gamma\setminus \{z\to y'\}}\rho_{z,r}(x)=1,
\end{align*}
where $\tau_z$ denotes the first hitting time to state $z$. This proves \eqref{dswr2} and thus completes the proof of Lemma \ref{lmm_dswr}. \qed
 
\subsection{Proof of Lemma \ref{lmm_scb}}

Due to Corollary \ref{lmm_odcb} and Lemma \ref{lmm_dswr}, it suffices to prove one direction. Suppose that $\pi$ is a complex-balanced distribution of $(\cN,\lambda)$ on $\Gamma$, then we need to verify  \eqref{def_vcbm} for $(\cN_1,\lambda_1)$. Let $\eta\in \cC\setminus \{z\}\subseteq \cC_1$, then from \eqref{def_lmd1out}, it follows that for any $x\in \Gamma$,
\begin{align}\label{equ_scb1}
\pi(x)\sum_{y'\colon \eta\to y'\in \cR_1}\lambda_{1,r}(x)=&\pi(x)\Big(\sum_{y'\colon\eta\to y'\in \cR_1^0}\lambda_{1,y\to y'}(x)+\sum_{i=1}^{p_z}\lambda_{1,\eta\to (z,i)}(x)\Big)\nonumber\\
=&\pi(x)\Big(\sum_{y'\colon\eta\to y'\in \cR, y'\neq z}\lambda_{\eta\to y'}(x)+\lambda_{\eta\to z}(x)\Big)=\pi(x)\sum_{y'\colon\eta\to y'\in \cR}\lambda_{y\to y'}(x).
\end{align} 
As $\pi$ is complex-balanced for $(\cR,\lambda)$ and $\phi_1=\phi\circ\psi_1=\phi$ on $\cC_1\setminus \{(z,i)\colon i=1,\dots, p_z\}$, we have 
\begin{align}\label{equ_scb2}
\pi(x)\sum_{y'\colon \eta\to y'\in \cR}\lambda_{y\to y'}(x)=&\sum_{y\colon y\to\eta\in \cR}\pi\big(x+\phi(y)-\phi(\eta)\big) \lambda_{y\to\eta}\big(x+\phi(y)-\phi(\eta)\big)\\
=&\sum_{y\colon y\to\eta\in \cR_1^0}\pi\big(x+\phi_1(y)-\phi_1(\eta)\big) \lambda_{y\to\eta}\big(x+\phi_1(y)-\phi_1(\eta)\big)\nonumber\\
&\quad+\pi\big(x+\phi_1(z)-\phi_1(\eta)\big) \lambda_{z\to \eta}\big(x+\phi_1(z)-\phi_1(\eta)\big).\nonumber
\end{align}
Then, \eqref{def_vcbm} is a consequence of \eqref{dswr1}, \eqref{equ_scb1} and \eqref{equ_scb2}. Next, we will show \eqref{def_vcbm} for $\eta=(z,i)$ with $i\in \{1,\dots, p_z\}$. Without loss of generality, assume $i=1$.

By definition,  $y_1\ce{->}(z,1)$ is the only incoming reaction of $(z,1)$. Therefore, if $x\not\geq \phi(z)$, then Condition \ref{def_cptsto} and the fact that $\phi_1=\phi\circ \psi_1$ with $\psi_1$ given by \eqref{def_vphi1} yields 
\[
0=\sum_{y'\colon (z,1)\to y'\in \cR_1^{out}}\pi(x)\lambda_{1,(z,1)\to y'}(x)=\pi\big(x+\phi_1(y_1)-\phi_1((z,1))\big)\lambda_{y_1\to (z,1)}\big(x+\phi_1(y_1)-\phi_1((z,1))\big).
\] 
Otherwise, assume $x\geq \phi(z)$. Let $\cX$ be the set of all complexes in $\cC$ that are in the same connected component 
 of $z$. As $(\cN,\lambda)$ is complex-balanced under $\pi$, then for any $\eta\in \cX$, 
\begin{align*}
0<\pi(x')\sum_{y'\colon \eta\to y'\in \cR}\lambda_{\eta\to y'}(x')=\sum_{y\colon y\to \eta\in \cR}\pi\big(x'+\phi(y)-\phi(\eta)\big)\lambda_{y\to \eta}\big(x'+\phi(y)-\phi(\eta)\big),
\end{align*}
where $x'=x+\phi(\eta)-\phi(z)\geq \phi(\eta )$. This yields that 
\begin{align*}
1=\frac{\sum_{y\colon y\to \eta\in \cR}\pi\big(x+\phi(y)-\phi(z)\big)\lambda_{y\to \eta}\big(x+\phi(y)-\phi(z)\big)}{\pi\big(x+\phi(\eta)-\phi(z)\big)\sum_{y'\colon \eta\to y'\in \cR}\lambda_{\eta\to y'}\big(x+\phi(\eta)-\phi(z)\big)}.
\end{align*}
We construct an irreducible DTMC  taking values in the finite state space $\cX$ with transition probabilities
\[ p_{z_1,z_2}=\mathbb{P}_{z_1}(z_2)=\frac{\pi\big(x+\phi(z_2)-\phi(z)\big)\lambda_{z_2\to z_1}\big(x+\phi(z_2)-\phi(z)\big)}{\pi\big(x+\phi(z_1)-\phi(z)\big)\sum_{y'\colon z_1\to y'\in \cR}\lambda_{z_1\to y'}\big(x+\phi(z_1)-\phi(z)\big)},
\]
for any $z_1,z_2\in\cX$. Then, $ p_{z_1,z_2}>0$ if and only if $z_2\ce{->}z_1\in \cR$. Thus, the chain is recurrent. With $\tau_z$ denoting the first hitting time to state $z$, we have $\mathbb{P}_{y_1}(\tau_z<\infty)=1$.

This proves \eqref{def_vcbm} with $\eta=(z,1)$, if it holds that
\begin{align}\label{equ_tauz}
\mathbb{P}_{y_1}(\tau_z<\infty)=\frac{\sum_{(z,1)\to y'\in \cR_1}\lambda_{1,(z,1)\to y'}(x)\pi(x)}{\pi\big(x+\phi_1(y_1)-\phi_1((z,1)\big)\lambda_{1,y_1\to (z,1)}\big(x+\phi_1(y_1)-\phi_1((z,1)\big)}.
\end{align}
First, by definition it is clear that
\[
\mathbb{P}_{y_1}(\tau_z<\infty)=p_{y_1,z}+\sum_{z'\in \cX\setminus\{z\}}p_{y_1,z'}p_{z',z}+\sum_{k=2}^{\infty}\sum_{\{z_1,\dots,z_k\}\subseteq \cX\setminus\{z\}}p_{y_1,z_1}\Big(\prod_{i=1}^{k-1}p_{z_{i},z_{i+1}}\Big)p_{z_{k},z}
\]
and
\[
R=\frac{\sum_{y'\in \cC}\sum_{\gamma \in \Gamma_{y_i\to z\to y'}}\prod_{r'\in \gamma\setminus \{z\to y'\}}\rho_{z,r'}(x)\lambda_{z\to y'}(x)}{\pi\big(x+\phi(y_1)-\phi(z)\big)\lambda_{y_1\to z}\big(x+\phi(y_1)-\phi(z)\big)},
\]
where $R$ denotes that right hand side of \eqref{equ_tauz}.
Additionally, for any $z_3\in \cX$, such that $z_1\ce{->}z_3\in \cR$, it follows from \eqref{def_rhosto} that
\begin{align*}
p_{z_1,z_2}=\frac{\pi\big(x+\phi(z_2)-\phi(z)\big)\lambda_{z_2\to z_1}\big(x+\phi(z_2)-\phi(z)\big)}{\pi\big(x+\phi(z_1)-\phi(z)\big)\lambda_{z_1\to z_3}\big(x+\phi(z_1)-\phi(z)\big)}\rho_{z,z_1\to z_3}(x).
\end{align*}
Consequently,
\begin{align}\label{equ_tauz1}
p_{y_1,z}=&\frac{\pi(x)\lambda_{z\to y_1}(x)\rho_{z,y_1\to z}(x)}{\pi\big(x+\phi(y_1)-\phi(z)\big)\lambda_{y_1\to z}\big(x+\phi(y_1)-\phi(z)\big)},
\end{align}
\begin{align*}
p_{y_1,z'}p_{z',z}=&\frac{\pi\big(x+\phi(z')-\phi(z)\big)\lambda_{z'\to y_1}\big(x+\phi(z')-\phi(z)\big)}{\pi\big(x+\phi(y_1)-\phi(z)\big)\lambda_{y_1\to z}\big(x+\phi(y_1)-\phi(z)\big)}\rho_{z,y_1\to z}(x)\nonumber\\
&\times\frac{\pi(x)\lambda_{z\to z'}(x)}{\pi\big(x+\phi(z')-\phi(z)\big)\lambda_{z'\to y_1}\big(x+\phi(z')-\phi(z)\big)}\rho_{z,z'\to y_1}(x)\nonumber\\
=&\frac{\pi(x)\lambda_{z\to z'}(x)\rho_{z,y_1\to z}(x)\rho_{z,z'\to y_1}(x)}{\pi\big(x+\phi(y_1)-\phi(z)\big)\lambda_{y_1\to z}\big(x+\phi(y_1)-\phi(z)\big)},
\end{align*}
and by iteration, letting $z_0=y_1$,
\begin{align}\label{equ_tauz3}
p_{y_1,z_1}\Big(\prod_{i=1}^{k-1}p_{z_{i},z_{i+1}}&\Big)p_{z_{k},z}
= \frac{\pi(x)\lambda_{z\to z_k}(x)\rho_{z,y_1\to z}(x)\prod_{i=1}^{k}\rho_{z,z_i\to z_{i-1}}(x)}{\pi\big(x+\phi(y_1)-\phi(z)\big)\lambda_{y_1\to z}\big(x+\phi(y_1)-\phi(z)\big)},
\end{align}
for all $k\geq 2$. Then,  \eqref{equ_tauz} follows from \eqref{equ_tauz1}-\eqref{equ_tauz3} and the definition of $\Gamma_{y_1\to z\to y'}$. The proof of this lemma is complete.  \qed

\subsection{Proof of Lemma  \ref{lmm_cfnt}}

 If $\cN_{M-1} = \cN_M$, then by definition, every complex in $\cC_{M-1}'\cap \cC' = \cC_{M}'\cap \cC'$ has only one incoming reaction. On the other hand, if $\cN_{M-1}\neq \cN_M$, then the $M$ complexes in $\cC'$ are cleaved sequentially in $\cN_1,\dots \cN_M$, and thus $\cC_M'\cap \cC' \subseteq \cC_M \cap \cC' = \emptyset$. Therefore, in either case, no complex in $\cC'\cap \cC_M$ has multiple incoming reactions. We will show that if $(y,i)\in \cC_M'$ is a copy of $y\in \cC'$, then $(y,i)$ has only one incoming reaction in $\cN_M$. First, $y\in \cC'$ has only one incoming reaction in $\cN$, but  multiple incoming reactions in $\cN_{m-1}$ for some $m\in\{1,\dots, M\}$; otherwise $(y,i)$ is not  in $\cC_M'\subseteq \cC_M$. Recall that when one-node cleaves a complex, only the incoming reactions of complexes that are products of the cleaved complex might change. It follows that the multiple incoming reactions in $\cR_{m-1}$ are due to the cleaving of a complex $y' \in \cC'\cup \{z\}$ in $\cN_{m'}$ with $m'<m$, that is, the reactant of the only incoming reaction of $y$ in $\cR,\dots,\cR_{m'-1}$. After cleaving $y'$, the reactant of each incoming reaction of $y$ in $\cC_{m'},\dots,\cC_{m-1}$ is a copy of $y'$, and thus when cleaving $y$ in $\cN_m$, the reactant $(y',j)$ of the only incoming reaction of $(y,i)$ in $\cR_m$ is a copy of $y'$. As in the cleaving iteration, only complexes in $\cC'$ might be cleaved. The copy $(y', j)$ is not  cleaved in $\cN_{m+1},\dots,\cN_M$. As a consequence, the incoming reactions of $(y,i)$ will not change, namely, $(y,i)$ has only one incoming reaction $(y',j)\ce{->}(y,i)$ in $\cR_{m+1},\dots, \cR_M$.

The only concern now is the cleaving of a complex $y$ in $\cN_m$ with some $m\in\{1,\dots, M\}$, fulfilling $y\ce{->}(z,i)\in \cR_{m-1}$ for some $i=1,\ldots,p_z$. The situation is illustrated in Figure \ref{fig:fig2}. Consider the RN $\cN$. Complex $z$ has two incoming reactions, and $\cC' = \{y_1,y_2\}$, in which each complex has only one incoming reaction. The only cycle including $y_2\ce{->}z$ included in $\cN$ is $\{z\ce{->}y_1\ce{->}y_2\ce{->}z\}$. Thus, after cleaving $z$ in $\cN_0$, the complex $(z,1)$ has only one outgoing reaction $(z,1)\ce{->}y_1$. Then, $y_1$ is cleaved in the same manner, resulting in the cleaved RN $\cN_1$. It remains to cleave the complex $y_2$. Note that there is only one cycle including $y_2\ce{->}(z,1)$ in $\cN_1$. Therefore, after cleaving $y_2$, the complex $(z,1)$ has only one incoming reaction in $\cR_2$. This observation allows us to complete the proof as follows.

Assume $y\in\cC_{m-1}\cap \cC'$ has multiple incoming reactions in $\cR_{m-1}$ (for some $m$), $y\ce{->}(z,1)\in \cR_{m-1}$, and $y$ is cleaved in $\cN_m$.

The reaction $y\ce{->}(z,1)$ is a result of the cleaving of $z$ in $\cN_0$, that is $y\ce{->}z\in\cR$ and $y\ce{->}(z,1)\in \cR_0$ is the only incoming reaction of $(z,1)$ in $\cR_0$. Additionally, $y\in \cC'$ has only one incoming reaction in $\cR$. It follows that the multiple incoming reactions of $y$ in $\cR_{m-1}$ come from the cleaving of some complex $y'$ in $\cN_{m'}$ with $m'\in\{0,\dots, m-1\}$. By iteration, we  find a sequence of reactions
\[
\{z\ce{->}y^{(1)}\ce{->}\cdots \ce{->}y^{(k)}\ce{->}y\}\subseteq \cR
\]
with $k\in\{0,\dots, m-1\}$, such that for each $i\in\{0,\dots, k\}$, complex $y^{(i)}\in \cC'$ is cleaved in $\cN_{m_i}$ that increases the number of incoming reactions of $y^{(i+1)}$in $\cR_{m_i}$, where $m_0,\dots, m_k$ are non-negative integers fulfilling $0=m_0<m_1<\dots<m_k\leq m-1$, with the convention that $y^{(0)}=z$. Because $\{y,y^{(1)},\dots, y^{(k)}\}\subseteq \cC'$, by definition $y^{(i)}\ce{->}y^{(i+1)}$ is the only incoming reaction of $y^{(i+1)}$ in $\cR$ for all $i=0,\dots,k$, where $y^{(0)}=z$ and $y^{(k+1)}=y$. Thus,
\[
\gamma_0=\{z\ce{->}y^{(1)}\ce{->}\cdots \ce{->}y^{(k)}\ce{->}y\ce{->}z\}
\]
is the only cycle including reaction $y\ce{->}z$ in $\cN$. Therefore, $(z,1)\ce{->}y^{(1)}$ is the only outgoing reaction of $(z,1)$ in $\cR_0$.

Next, consider the cleaving of $y^{(1)}$ in $\cN_{m_1}$. Neither $(z,1)$ nor $y^{(1)}$ is cleaved in $\cN_1,\dots, \cN_{m_1-1}$, thus $(z,1)\ce{->}y^{(1)}\in \cR_{m_1-1}$. Therefore, after the cleaving of $y^{(1)}$ in $\cN_{m_1}$, there is a copy of $y^{(1)}$, denoted by $(y^{(1)},1)$ in $\cC_{m_1}$ such that $(z,1)\ce{->}(y^{(1)},1)\in \cR_{m_1}$. Additionally, this reaction is the only incoming reaction of $(y^{(1)},1)$ and the only outgoing reaction of $(z,1)$ in $\cR_{m_1}$. Similarly, since neither $(y^{(1)},1)$ or $(z,1)$ is cleaved in $\cN_{m_1+1},\dots, \cN_m$, then it follows that $(z,1)\ce{->}(y^{(1)},1)\in \cR_{m_1}$ is the only incoming reaction of $(y^{(1)},1)$ and the only outgoing reaction of $(z,1)$ in $\cR_{m_1+1},\dots, \cR_{m-1}$ as well. On the other hand, because none of $ y^{(1)},\dots, y^{(k)}, (z,1)$ are cleaved in $\cN_1,\dots, \cN_{m_1-1}$, it holds that
\[
\gamma_1=\{y^{(1)}\ce{->}y^{(2)}\ce{->}\cdots \ce{->}y^{(k)}\ce{->}y\ce{->}(z,1)\ce{->}y^{(1)}\}.
\]
is also the only cycle including $(z,1)\ce{->}y^{(1)}$ in $\cN_{m_1-1}$. Thus, $(y^{(1)},1)\ce{->}y^{(2)}$ is the only outgoing reaction of $(y^{(1)},1)$ in $\cR_{m_1}$, and thus in $\cR_{m_1+1},\dots \cR_{m_2-1}$. 

By iteration, after cleaving $y^{(k)}$, there is a sequence 
\[
\{(z,1)\ce{->}(y^{(1)},1)\ce{->}\cdots \ce{->}(y^{(k)},1)\ce{->}y\}\subseteq \cR_{m_k}
\]
such that $(y^{(i)},1)\ce{->}(y^{(i+1)},1)$ is the only outgoing reaction of $y^{(i)}$ for all $i=0,\dots, k$ in $\cR_{m_k},\dots, \cR_{m-1}$, where $y^{(0)}=z$ and $(y^{(k+1)},0)=y$. This implies that the only cycle including $y\ce{->}(z,1)$ in $\cN_{m-1}$ is
\[
\gamma_2=\{y\ce{->}(z,1)\ce{->}(y^{(1)},1)\ce{->}\cdots \ce{->}(y^{(k)},1)\ce{->}y\}.
\]
As a consequence, complex $(z,1)$ has only one incoming reaction in $\cR_{m}$, provided $y\ce{->}(z,1)$ is the only incoming reaction of $(z,1)$ in $\cR_{m-1}$. 
This proves that the number of incoming reactions of $(z,1)$ is one in $\cR_m$. The proof of this lemma is thus complete. \qed

\section{Conflict of interest disclosure}

The authors declare they have no conflict of interest relating to the content of this article.

\section{Funding}
The work presented in this article is supported by Novo Nordisk Foundation (Denmark), 
grant NNF19OC0058354.


\end{document}